\documentclass[11pt]{amsart}

\usepackage{amsmath}
\usepackage{amsfonts}
\usepackage{amsthm}
\usepackage{euscript}
\usepackage{amssymb}
\usepackage{graphicx}
\usepackage{mathrsfs}
\usepackage{enumitem}
\usepackage{indentfirst}
\usepackage{csquotes}
\usepackage{bbm}
\usepackage{xcolor}

\newcommand{\field}[1]{\mathbb{#1}}
\newcommand{\N}{\field{N}}

\newcommand{\la}{\lambda}

\newcommand{\tb}{\widetilde \beta}

\renewcommand{\epsilon}{\varepsilon}

\numberwithin{equation}{section}
\newtheorem{theorem}{Theorem}[section]
\newtheorem{lemma}{Lemma}[section]
\newtheorem{corr}[theorem]{Corollary}

\newtheorem{proposition}{Proposition}[section]
\newtheorem{deff}[theorem]{Definition}
\newtheorem{remark}[theorem]{Remark}
\newcommand{\bth}{\begin{theorem}}
\newcommand{\ble}{\begin{lemma}}
\newcommand{\bcor}{\begin{corr}}

\newcommand{\bdeff}{\begin{deff}}

\newcommand{\bprop}{\begin{proposition}}
\newcommand{\ele}{\end{lemma}}
\newcommand{\ecor}{\end{corr}}
\newcommand{\edeff}{\end{deff}}

\newcommand{\eprop}{\end{proposition}}

\newcommand{\cd}{\, \cdot\, }

\newcommand{\e}{\varepsilon}

\newcommand{\supp}{\text{supp }}
\renewcommand{\Pi}{\varPi}

\renewcommand{\epsilon}{\varepsilon}

\newcommand{\R}{{\mathbb R}}

\newcommand{\1}{\mathbbm{1}}

\newcommand{\Atn}{A^{\theta_0}_\nu}

\newcommand\norm[1]{\left\lVert#1\right\rVert}

\newcommand{\lpnorm}[3]{\norm{#1}_{L^{#2}(#3)}}
\newcommand{\altlpnorm}[2]{\norm{#1}_{#2}}
\newcommand{\mixlpnorm}[4]{\norm{#1}_{L^{#2}_{t} L^{#3}_{x}(#4)}}

\newcommand{\mixlitbiglpnorm}[5]{\norm{#1}_{\ell ^{#2}_{#3} L^{#4}(#5)}}
\newcommand{\othermixlitbiglpnorm}[8]{\norm{#1}_{L_{#2}^{#3}\ell ^{#4}_{#5} L_{#6}^{#7}({#8})}}

\newcommand{\Lpopnorm}[5]{\norm{#1}_{L^{#2}(#3) \to L^{#4}(#5) }}

\newcommand{\Z}[1][]{\mathbb{Z}^#1}

\newcommand{\AvtO}{A_{\nu}^{\theta_0}}

\newcommand{\laplac}[1]{\Delta_{#1}}
\newcommand{\Sl}{S_{\la}}

\newcommand{\scaledsol}[1]{e^{-i\la^{-1}#1\laplac{g}}}
\newcommand{\scaledsoll}[1]{e^{i\la^{-1}#1\laplac{g}}}
\newcommand{\bilinKT}[2]{\mathcal{B}_{\la #1}(#2, \{\phi_{\nu}\})}

\begin{document}

\title[Inhomogeneous Strichartz Estimates]{Inhomogeneous Strichartz Estimates on Manifolds with nonpositive curvature and applications}

\author{Xiaoqi Huang}
\address[X.H.]{Department of Mathematics, Louisiana State University, Baton Rouge, LA 70808}
\email{xhuang49@lsu.edu}

\author{Connor Quinn}
\address[C.Q.]{Department of Mathematics,  Johns Hopkins University,
Baltimore, MD 21218}
\email{cquinn24@jhu.edu}

\begin{abstract}
    The main purpose of this paper is to establish lossless inhomogeneous Strichartz estimates for solutions to the Schrödinger equation on compact manifolds with nonpositive curvature over frequency-dependent short time intervals. As applications, we improve upon the Sobolev norm growth bounds for the cubic NLS on 3-dimensional compact manifolds established in \cite{planchon2017growth}, and extend the homogeneous logarithmic Strichartz estimates of \cite{huang2024strichartz} to Schrödinger operators with critically singular potentials.
\end{abstract}

\maketitle

\section{Introduction}

Let $(M,g)$ be a compact Riemannian manifold of
dimension $d\ge2$, $\Delta_g$ be
the associated Laplace-Beltrami operator and let
\begin{equation}\label{i1}
u(x,t)=e^{it\Delta_g}f(x)
\end{equation}
be the solution to the Schr\"odinger equation
on $M\times \R$,
\begin{equation}\label{i2}
i\partial_t u(x,t)+\Delta_g u(x,t)=0 \ , \quad u(x,0)=f(x).
\end{equation}
For an interval $I \subset \R$, define 
$$\mixlpnorm{u}{p}{q}{M \times I}=\Bigl(\, \int_I \, \|u(\, \cdot\, , t)\|_{L^q_x(M)}^p \, dt\,\Bigr)^{1/p}.$$
Recall that in the work of Burq, G\'erard and Tzvetkov~\cite{bgtmanifold}, the authors proved mixed-norm Strichartz estimates
\begin{equation}\label{i3}
\|u\|_{L^p_tL^q_x(M\times [0,1])}
\lesssim \|f\|_{H^{1/p}(M)}
\end{equation}
for all {\em admissible} exponent pairs $(p,q)$.  By this
we mean, as in Keel and Tao~\cite{KT}, that
\begin{equation}\label{i4}
d\big(\tfrac12-\tfrac1q\big)=\tfrac2p \,  \, \, \text{and}
 \, \, \, 2< q\le \tfrac{2d}{d-2} \, \,
\text{if } \, d\ge 3, \, \, \, 
\text{or } \, 2< q<\infty \, \, 
\text{if } \, \, d=2.
\end{equation}
Also, in \eqref{i3} 
$H^\mu$ denotes the
standard Sobolev space
\begin{equation}\label{i5}
\|f\|_{H^\mu(M)}=
\bigl\| \, (I+P)^\mu f\, \bigr\|_{L^2(M)}, 
\quad \text{with } \, \, P=\sqrt{-\Delta_g},
\end{equation}
and ``$\lesssim$'' in \eqref{i3} and, in what
follows, denotes an inequality with an implicit,
but unstated, constant $C$ which can change at each occurrence.

Due to the absence of a global dispersive estimate,
the proof of \eqref{i3} relies on a lossless estimate over frequency-dependent short time intervals. To state this, let us fix a Littlewood--Paley bump function $\beta$ satisfying
\begin{equation}\label{i6}
\beta \in C^\infty_0((1/2,2)), \,\beta=1 \,\,\text{on} \,\,[3/4, 5/4]
\quad \text{and} \quad
1=\sum_{k=-\infty}^\infty \beta(2^{-k}s),
\qquad s>0.
\end{equation}
Define $$\beta_0(s)=1-\sum_{k=1}^\infty
\beta(2^{-k}s)\in C^\infty_0(\R_+)$$ and
$\beta_k(s)=\beta(2^{-k}s)$, $k=1,2,\dots$. 
Then, by the standard Littlewood--Paley theory on compact manifolds (see, e.g., \cite{SFIO2}), together with arguments from the earlier works \cite{bgtmanifold, huang2024strichartz}, one can reduce \eqref{i3} to the following uniform estimate: 
\begin{equation}\label{i7}
\| e^{it\Delta_g}\beta(P/\la)f\|_{L^p_tL^q_x(M\times [0,1])}\lesssim \la^{\frac1{p}}\, \|f\|_{L^2(M)}, \quad \la \gg 1.
\end{equation}
Here and in what follows $\la\gg 1$ means $\la \ge C$ for some large constant.
Burq, G\'erard and Tzvetkov proved this estimate in
\cite{bgtmanifold} by showing that one always has the 
following uniform dyadic estimates over very small intervals:
\begin{equation}\label{i8}
\|e^{it\Delta_g} \beta(P/\la)f\|_{L^p_tL^q_x(M\times [0,\la^{-1}])}\lesssim\, \|f\|_{L^2(M)}, \quad \la \gg 1.
\end{equation}
It is not hard to see that 
 \eqref{i8}  yields 
\eqref{i7}, since one can write $[0,1]$ as the
union of $\approx \la$ intervals of length $\la^{-1}$
and thus obtain \eqref{i7} by adding up the
uniform estimates on each of these subintervals. Also, the  bounds in \eqref{i8} cannot be improved
on {\em any} manifold, as one can see by taking $f(x)=f_\la(x)=\beta(P/\la)(x,x_0)$ where $\beta(P/\la)(x,y)$ is the kernel of the 
Littlewood-Paley operators and $x_0$ is an arbitrary fixed point on $M$.  See also \cite{HSSchro} for another class of examples, concentrating near geodesics, which saturate \eqref{i8}.

However, on manifolds satisfying additional geometric assumptions, it is possible to extend the lossless estimate to time intervals longer than $\lambda^{-1}$.
For instance, if the manifold has nonpositive sectional curvature, Huang and Sogge \cite{huang2024strichartz} proved the following uniform bound:
\begin{equation}\label{i9}
\| e^{it\Delta_g} \beta(P/\la) f\|_{L^p_tL^q_x(M
\times [0, \, \log\la\cdot \la^{-1}])}\lesssim \, \|f\|_{L^2(M)}, \quad \la\gg 1.
\end{equation}
Then \eqref{i9} yields an improvement over \eqref{i3} with a logarithmic gain in derivatives, which follows by writing $[0,1]$ as a union of $\approx \la / \log \la$ intervals and adding up the uniform estimate \eqref{i9} on each of these intervals.
The main ideas in the work \cite{huang2024strichartz} are based on earlier work by 
 Blair, Huang and Sogge \cite{blair2023strichartz}, where a similar result with a smaller logarithmic gain in derivatives was obtained.

Under more specialized geometries, such as the flat torus, it is possible to extend \eqref{i9} to much larger time intervals of length $\lambda^{-\delta}$ for some fixed $0<\delta<1$; see the recent work of the second author \cite{quinn2026}. For the special case $p=q=\frac{2(d+2)}{d}$, when $d=2$, Herr and Kwak \cite{herr2023strichartz} obtained a lossless version of \eqref{i9} on time intervals of length $(\log \lambda)^{-1}$; see also Skouloudis and Yu \cite{skouloudis2026strichartz} for a related result on intervals of length $(\log \lambda)^{-C}$ for $d=1$. In higher dimensions, the $\ell^2$-decoupling theorem of Bourgain and Demeter \cite{BourgainDemeterDecouple} yields an almost lossless version of \eqref{i9} on unit time intervals with only a $\lambda^\varepsilon$ loss of derivatives; see also Deng, Germain and Guth \cite{DGG} for further discussion on irrational tori.

Now consider the inhomogeneous Schr\"odinger equation
\begin{equation}
    i\partial_{t}u + \laplac{g}u = F, \ \ \  u(x,0) =0.
\end{equation}
By Duhamel's principle, the solution is given by
\begin{equation} \label{inhom sol}
    u(x,t) =  \int_{0}^{t}e^{i(t-s)\laplac{g}}F(x,s)ds.
\end{equation}
The first main result of this paper is the following inhomogeneous analogue of \eqref{i9}.
\begin{theorem}\label{thm}  Let $M$ be a compact $d$-dimensional manifold, $d\ge 2$, whose sectional curvatures are all nonpositive, and let $\beta$ be defined as in \eqref{i6}. Then for all admissible pairs $(p_1,q_1)$ and $(p_2,q_2)$ satisfying \eqref{i4} and $\la\gg 1$,
\begin{equation}\label{i10}
   \mixlpnorm{\int_0^t e^{i(t-s)\Delta_g}\beta(P/\la) F(\cdot, s)ds}{p_1}{q_1}{M \times [0,\log\la \cdot\la^{-1}]}  \lesssim  \|F\|_{L^{p'_2}_tL^{q_2'}_x(M\times [0, \log\la\cdot \la^{-1}])} 
\end{equation}
\end{theorem}
\noindent Here $p'$ denotes the conjugate exponent such that $\frac1p+\frac1{p'}=1$.
If the interval is replaced by $[0,\la^{-1}]$, then the above estimate follows from the dispersive estimate established in \cite{bgtmanifold} (see \eqref{s1b} below) along with the abstract theorem of Keel and Tao \cite{KT}.  Moreover, if $(p_1,p_2)\neq (2,2)$, then \eqref{i10} follows from the homogeneous estimate \eqref{i9} combined with the Christ--Kiselev lemma. Therefore, the main new contribution of Theorem~\ref{thm} concerns the double-endpoint case where $$d\ge 3,\,\,\,\text{and}\,\,\, (p_1,q_1)=(p_2,q_2)=(2, \tfrac{2d}{d-2}).$$ 

The main ideas in the proof of Theorem~\ref{thm} are similar to those in \cite{huang2024strichartz}. More precisely, we iteratively apply a height decomposition twice, and in the low-height regions arising at each stage, we use an almost orthogonality inequality established in \cite{huang2024strichartz}, involving microlocal pseudodifferential cutoff operators in both space and frequency. This allows us to reduce the inhomogeneous Strichartz estimate to its analog, which involves  microlocal pseudodifferential operators. We then prove the microlocalized inhomogeneous estimate using the microlocalized dispersive estimate established in \cite{blair2023strichartz} and
adapting the proof of the main theorem in \cite{KT}.

It is possible to extend Theorem~\ref{thm} to noncompact manifolds with bounded geometry and nonpositive curvature. For the homogeneous estimate \eqref{i9}, this was carried out in \cite{huang2026lossless}, where the authors further used this estimate as one of the main ingredients to upgrade the lossless estimate on frequency-dependent intervals to a global-in-time estimate on asymptotically hyperbolic surfaces with negative curvature. We hope to further explore applications of Theorem~\ref{thm} in the noncompact setting in a future work.

Strichartz estimates have a broad range of applications to nonlinear 
Schr\"odinger (NLS) equations; see, e.g., Cazenave \cite{cazenave2003semilinear} 
for a comprehensive treatment. On compact manifolds, the estimates of 
\cite{bgtmanifold} have been used to study well-posedness for subcritical 
NLS. In a related direction, the lossless inhomogeneous Strichartz estimate 
over short time intervals of size $\lambda^{-1}$ from \eqref{i10} plays a 
key role in the work of Planchon, Tzvetkov and Visciglia 
\cite{planchon2017growth} on the growth of $H^{s}$ $(s > 1)$ norms for the cubic NLS on 3-dimensional compact manifolds. Such estimates have received much interest over the last decade due to their connection with a phenomenon known as ``weak-wave turbulence", by which solutions to NLS transfer energy from low to high frequencies over time. This causes the $H^{s}$ norm of the solutions to grow, while the energy remains bounded. See the introductions to the papers \cite{deng2018growth} and \cite{planchon2017growth} and the references therein for more details. On manifolds with nonpositive sectional curvature, 
Theorem~\ref{thm} provides a lossless estimate over the longer 
interval $[0, \log\lambda \cdot \lambda^{-1}]$, and we can utilize this to attain 
improved bounds on the growth of $\|u(\cdot, t)\|_{H^m}$ in this setting over the results in \cite{planchon2017growth}. To be more specific, consider the cubic NLS
\begin{equation} \label{cubic NLS}
    i\partial_t u + \laplac{g}u = |u|^2u \ , \ \ \ u(x,0) = f(x).
\end{equation} 
In \cite{planchon2017growth}, the authors showed that for any $m \geq 2$ the $H^{m}(M)$ norm of $u(\cdot,t)$ is bounded by $\exp(CT)$ for any $t \in [0,T]$. Using the double-endpoint version of Theorem~\ref{thm} and adapting the method used in \cite{planchon2017growth}, we are able to prove the following theorem.
\begin{theorem}\label{thm4}
    Let $(M,g)$ be a $3$-dimensional compact manifold with nonpositive sectional curvature. Suppose $u \in C_{t}(H^{m}(M))$ is a solution to \eqref{cubic NLS}. If $m \in \N$ with $m \geq 2$ and $f \in H^{m}(M)$ then for any $T > 0$ we have
    \begin{equation}
        \sup_{0 \leq t\leq T}\lVert u(\cdot,t)\rVert_{H^{m}(M)} \leq C\exp{\big(C\sqrt{T}\big)}
    \end{equation}
   where $C = C(m,\lVert f\rVert_{H^{m}}).$
\end{theorem}
\noindent 

As another immediate application of Theorem~\ref{thm}, we are also able to extend the homogeneous estimate \eqref{i9} to Schr\"odinger equations with critically singular potential terms.
Let $V \in L^{d/2}(M)$ be real-valued. Set $H_V = -\laplac{g} + V(x)$ and consider the Schr\"{o}dinger equation
\begin{equation}
    i\partial_{t}u - H_{V}u = 0, \ \ u(x,0) = f(x).
\end{equation}
Then we can write
\begin{equation}
    u(x,t) = e^{-itH_V}f(x).
\end{equation}
By adding a constant to $V$ if necessary, we may assume $H_V$ is positive and denote $P_V=\sqrt{H_V}$.
\begin{theorem}\label{thm2}  Let $V$ be as above, $\beta$ be as in \eqref{i6}, and let $(M,g)$ be a compact manifold with nonpositive sectional curvature. Assume that the pair of exponents $(p,q)$ satisfies \eqref{i4}. Then
\begin{equation}\label{i11}
\| e^{-itH_V} \beta(P_V/\la) f\|_{L^p_tL^q_x(M
\times [0, \, \log\la\cdot \la^{-1}])}\lesssim \, \|f\|_{L^2(M)}, \quad \la\gg 1.
\end{equation}
\end{theorem}
Here  we call $V \in L^{d/2}(M)$ critically singular
potentials, since multiplication by elements of $L^{d/2}$ scale the same as the Euclidean
Laplacian.
The proof of \eqref{i11} adapts the arguments in \cite{HSSchro}, and uses the double-endpoint version of Theorem~\ref{thm} together with a perturbative argument in which the additional term involving the potential is treated as an inhomogeneous forcing term. In particular, it is crucial to the argument, especially in the critical case, that \eqref{i10} holds without any loss of derivatives. As in earlier works such as \cite{HSSchro,IvanoviciMT}, we also make use of the fact that the spatial frequency cutoff $\beta(\sqrt{H_V}/\lambda)$ can be interchanged with a suitable dyadic time-frequency cutoff. This allows us to view $ e^{itH_V} \beta(\sqrt{H_V}/\la) f$ as a frequency-localized solution to an inhomogeneous Schr\"{o}dinger equation associated with the standard Laplacian.

By using the Littlewood--Paley theory for the Schr\"{o}dinger operator $H_V$ (see \cite[Corollary 5.2]{HSSchro}) and summing the estimates over disjoint intervals of length $\lambda^{-1}\log\lambda$ covering $[0,1]$, we obtain the following consequence of \eqref{i11}.
\begin{theorem}\label{thm3} Under the assumptions of Theorem~\ref{thm2}, we have 
\begin{equation}
     \mixlpnorm{e^{-itH_V}f}{p}{q}{M \times [0,1]} \lesssim \lpnorm{(I + P_V)^{\frac{1}{p}}\log (2I + P_V)^{-\frac1p}f}{2}{M}.
\end{equation}
\end{theorem}

This paper is organized as follows. Sections 2-5 focus on the proof of Theorem~\ref{thm}. We begin Section 2 by making some initial reductions  to reduce ~\eqref{inhom sol} to a form that is more amenable to the use of the pseudodifferential cutoff operators we will ultimately use. Here we also define the large and small height regions mentioned earlier. Section 3 focuses on the proof of our desired estimate on the small height region, while Section 4 focuses on the proof on the large height region. In Section 5 we prove the microlocalized version of the double-endpoint inhomogeneous estimate mentioned above, which concludes the proof of Theorem~\ref{thm}. In Section 6 we sketch the proof of Theorem~\ref{thm4}, mainly focusing on the alterations to the argument in \cite{planchon2017growth} that yield our improved result on manifolds with nonpositive curvature. Finally, in Section 7, we prove Theorem~\ref{thm2}.

\section*{Acknowledgments} 
The authors would like to thank Professor Nikolay Tzvetkov for suggesting in a communication with Chris Sogge that was shared with us how combining the logarithmic improvement for the Strichartz estimates shown by Blair, Huang and Sogge in \cite{blair2023strichartz} with the techniques in the paper \cite{planchon2017growth} by Planchon, Tzvetkov and Visciglia could lead to improved Sobolev growth estimates for the cubic NLS equation on compact manifolds with nonpositive curvature.

\section{Initial Reductions for the proof of Theorem~\ref{thm}.}
As remarked earlier, the inhomogeneous estimate follows from the homogeneous one via the Christ--Kiselev lemma, except at the endpoint. Let $q_e=\frac{2d}{d-2}$.
 By a simple rescaling argument, the endpoint estimate in \eqref{i10} is equivalent to the following estimate
\begin{equation} \label{a}
    \mixlpnorm{\int_{0}^{t}\scaledsoll{(t-s)}\beta(P/\la)F(\cdot,s)ds}{2}{q_e}{M \times [0,\log\la]} \lesssim \la \mixlpnorm{F}{2}{q_e'}{M \times [0,\log\la]}.
\end{equation}

Let $T = c_0\log\la$ for some sufficiently small $c_0$ to be specified later. 
We first claim that \eqref{a} follows from
\begin{equation} \label{s1}
    \mixlpnorm{\int_{0}^{t}\scaledsoll{(t-s)}\beta(P/\la)F(\cdot,s)ds}{2}{q_e}{M \times [0, T]} \lesssim \la \mixlpnorm{F}{2}{q_e'}{M \times [0, T]}.
\end{equation}
To verify the claim, we divide the interval $[0,\log\la]$ into finitely many subintervals $[jT,(j+1)T],$ $j\in \mathbb{Z}, |j|\lesssim 1$. On each $I_j$, we estimate 
\begin{equation} \label{a1}
\begin{aligned}
     &\mixlpnorm{\int_{0}^{t}\scaledsoll{(t-s)}\beta(P/\la)F(\cdot,s)ds}{2}{q_e}{M \times [jT,(j+1)T]} \\
     &\le  \mixlpnorm{\int_{jT}^{t}\scaledsoll{(t-s)}\beta(P/\la)F(\cdot,s)ds}{2}{q_e}{M \times [jT,(j+1)T]} \\
     &\quad +\mixlpnorm{\int_{0}^{jT}\scaledsoll{(t-s)}\beta(P/\la)F(\cdot,s)ds}{2}{q_e}{M \times [jT,(j+1)T]}
\end{aligned}
\end{equation}
The first term on the right side can be handled using \eqref{a1} and a simple change of variables argument, while the second term follows from \eqref{i9} and duality. Summing over finitely many $j$ give us \eqref{a}.

Recall that the result of Burq-G\'erard-Tzvetkov from \cite{bgtmanifold} yields
\begin{equation}\label{s1a}
  \mixlpnorm{ \int_0^t \scaledsoll{(t-s)}\beta(P/\la) F(\cdot, s)ds}{2}{q_e}{M\times[0,1]} \lesssim \la \|F\|_{L^2_tL^{q_e'}_x([0,1]\times M)}.
\end{equation}
This is a consequence of the Keel-Tao theorem from \cite{KT}, the trivial $L^2 \to L^2$ estimate for the Schr\"{o}dinger propagator and the following dispersive type estimate, which was also originally established in \cite{bgtmanifold}:
\begin{equation}\label{s1b}
   \| \scaledsoll{(t-s)}\beta^2(P/\la) f\|_{L^{\infty}(M)} \lesssim \la^{\frac{d}{2}}|t-s|^{-\frac{d}{2}}\|f\|_{L^{1}(M)}, \ \ \ \  \text{if} \ \   |t-s|\lesssim 1. 
\end{equation}
See also \cite[Proposition 4.1]{blair2023strichartz} for a different proof of \eqref{s1b}.

We now turn to the proof of \eqref{s1}. Note that it suffices to consider the case where $$\|F\|_{L_{t}^2L_{x}^{q_e'}(M \times [0,T])} = 1.$$ Let $\1_+(t) $ denote the indicator function of $[0, +\infty)$. We can write 
\begin{equation}\label{defa}
    \int_0^t \scaledsoll{(t-s)}\beta(P/\la) F(\cdot, s)ds = \int_0^T \scaledsoll{(t-s)}\1_{+}(t-s)\beta(P/\la) F(\cdot,s)ds
\end{equation}
Now, let $A=\{(t,s)\in [0, T]^2: s \le [t]\}$, where $[t]$ denotes the integer part of $t$.  Let $\1_A$ be the indicator function of the set $A$. Then it is not hard to see that $\1_{+}(t-s)- \1_A(t,s)$ is supported in the region $|t-s|\le 1$. 
Thus, by \eqref{s1a} along with a change of variable argument, we have 
\begin{equation} \label{strichartz estimate for unit time error term}
    \begin{aligned}
        & \mixlpnorm{\int_0^T \scaledsoll{(t-s)}(\1_{+}(t-s)- \1_A(t,s))\beta(P/\la) F(\cdot, s)ds}{2}{q_e}{M \times [0,T]} \\
        &\le \bigg(\sum_{0\le j\le [T]}\mixlpnorm{ \int_0^T \scaledsoll{(t-s)}(\1_{+}(t-s)- \1_A(t,s))\beta(P/\la) F(\cdot, s)ds}{2}{q_e}{M \times [j,j+1]}^2 \bigg)^{\frac12
        } \\
        &= \bigg(\sum_{0\le j\le [T]}\mixlpnorm{\int_j^t \scaledsoll{(t-s)}\beta(P/\la) F(\cdot,s)\1_{[j,j+1]}(s)ds}{2}{q_e}{M \times [j,j+1]}^2\bigg)^{\frac12
        }\\
        &\lesssim \la \Big(\sum_{0\le j\le [T]}\mixlpnorm{F}{2}{q_e'}{M \times [j,j+1]}^2  \Big)^{\frac12 
        }\lesssim \la \mixlpnorm{F}{2}{q_e'}{M \times [0,T]}. 
    \end{aligned}
\end{equation}

Define the time-dilated integrated Schr\"{o}dinger operator
\begin{equation}
    S_\la F(x,t) =\int_0^T \scaledsoll{(t-s)} \1_A(t,s)\beta(P/\la) F(x,s)ds.
\end{equation}
Note that by \eqref{strichartz estimate for unit time error term} it suffices to establish
\begin{equation}\label{s4}
   \mixlpnorm{\Sl F}{2}{q_e}{M \times [0,T]} 
   \lesssim \la \mixlpnorm{F}{2}{q_e'}{M \times [0,T]}.
\end{equation}
We need to introduce several auxiliary operators that will allow us to apply  Proposition 2.3 in \cite{huang2024strichartz}, which is essential for our proof. The ``local'' auxiliary operators will be  the following ``quasimode'' operators adapted to the
scaled Schr\"odinger operators $i\la\partial_t +\Delta_g$, 
\begin{equation}\label{22.5}
\sigma_\la = \sigma\bigl(\la^{1/2}|D_t|^{1/2}-P\bigr) \, \tilde \beta(-D_t/\la),
\end{equation}
where
\begin{equation}\label{22.6}
\sigma\in {\mathcal S}(\R) \, \, \text{satisfies } \, \,
\sigma(0)=1 \, \, \text{and } \, \, \text{supp }\Hat \sigma
\subset \delta\cdot [1-\delta_0,1+\delta_0]=[\delta-\delta_0\delta, \delta+\delta_0\delta],
\end{equation}
with $0<\delta,\delta_0<1/8$ sufficiently small and
\begin{equation}\label{22.7}
\tilde \beta\in C^\infty_0((1/8,8)) \quad \text{such that } \, \,
\tilde \beta=1 \, \, 
\text{on } \, \, [1/6,6].
\end{equation}
Here, the properties of $\sigma$, as well as the small constants $\delta$ and $\delta_0$, are the same as those in \cite{blair2023strichartz}, except that we use $\beta(-D_t/\lambda)$ instead of $\beta(D_t/\lambda)$ in order to be consistent with the forward Schr\"odinger propagator $e^{it\Delta_g}$.
  
Note that we have
\begin{equation}\label{22.8}
I=\sum_{k=1}^{N_0} B_k(x,D),
\end{equation}
where each $B_k\in S^0_{1,0}(M)$ is a standard
pseudodifferential operator with symbol supported in
a small conic neighborhood of some $(x_k,\xi_k)\in
S^*M$. These operators and $N_0$ will not depend
on our parameter $\la\gg 1$.  Next, if 
$\tilde \beta$ is as in \eqref{22.7} then the
dyadic operators
\begin{equation}\label{22.9}
B=B_{k,\la} = B_k\circ \tilde \beta(P/\la)
\end{equation}
are uniformly bounded on $L^p$, i.e.
\begin{equation}\label{22.10}
\Lpopnorm{B}{p}{M}{p}{M}
=O(1) \quad \text{for } \, \, 1\le p\le \infty.
\end{equation}
The operators \(B_{k,\lambda}\) allow us to localize to a fixed coordinate chart on \(M\) and to use explicit formulae of parametrices for local operators such as \(\sigma_\lambda\) in our analysis.

Since $\tilde \beta(P/\la)S_\la F=S_\la F$, to prove \eqref{s4}, it suffices to show that for a fixed $B = B_{k,\la}$ we have
\begin{equation}\label{s5b}
   \mixlpnorm{B S_\la F}{2}{q_e}{M \times [0,T]} 
   \lesssim \la \mixlpnorm{F}{2}{q_e'}{M \times [0,T]}. 
\end{equation}

As in previous works of Blair-Huang-Sogge \cite{blair2023strichartz}, Huang-Sogge \cite{huang2024strichartz} and Quinn \cite{quinn2026} we decompose $M \times [0,T]$ into two subsets based on the height of $B\Sl F$.
Define
\begin{equation} \label{small height set def}
  A^{-}=\{(x,t)\in M \times [0,T]: |B S_\la F(x,t)|\lesssim \la^{\frac{d}{4}+\frac12+\e_0} \}
\end{equation}
and 
\begin{equation} \label{large height set def}
  A^{+}=\{(x,t)\in M \times [0,T]: |B S_\la F(x,t)|\gtrsim \la^{\frac{d}{4}+\frac12+\e_0} \}
\end{equation}
Here $\la^{\frac{d}{4}+\e_0}$ is the natural height cut off for the homogeneous case as in the earlier works \cite{blair2023strichartz, huang2024strichartz}. The additional $\la^{\frac12}$ factor can be viewed as arising from the dual homogeneous estimate
 $L^2_tL^{q_e'}_x \to L^2$. Indeed, the inhomogeneous estimate can essentially be regarded as the composition of the homogeneous estimate
$L^2 \to L^2_tL^{q_e}_x$ estimate with its dual version $L^2_tL^{q_e'}_x \to L^2$.

Clearly, \eqref{s5b} follows from the following two estimates:
\begin{equation} \label{small height sb}
    \mixlpnorm{BS_\la F}{2}{q_e}{A^{-}} 
   \lesssim \la \mixlpnorm{F}{2}{q_e'}{M \times [0,T]}
\end{equation}
and 
\begin{equation} \label{large height sb}
    \mixlpnorm{BS_\la F}{2}{q_e}{A^{+}} 
   \lesssim \la \mixlpnorm{F}{2}{q_e'}{M \times [0,T]}
\end{equation}
We will first focus on proving \eqref{small height sb}. The estimate \eqref{large height sb} is more involved than in previous works and we will come back to it after proving \eqref{small height sb}.

\section{Small Height Estimate}
To start, define $I_j = [j,j+1)\cap [0, T]$ for $0\le j\le [T]$. In addition, define $A^{-}_{j}=(M \times I_j) \cap A^{-}$ for $0 \leq j \leq [T]$. We then have
\begin{equation}
    \mixlpnorm{B S_\la F}{2}{q_e}{A^{-}} 
    = \Big( \sum_{j}\mixlpnorm{B S_\la F}{2}{q_e}{A^{-}_{j}}^2\Big)^\frac12.
\end{equation}
Note that by the definition of the set $A$, for $t\in I_j$, the $s$ support of $\1_A(t,s)$ is $[0, j]$. Thus, for $(x,t)\in A^{-}_j$ we have
\begin{equation}\label{s5a}
\begin{aligned}
    S_\la F(x,t)&= \int_0^{j} \scaledsoll{(t-s)} \beta(P/\la) F(x,s)ds     \\
    &= \scaledsoll{t}\left(\int_0^{j} e^{-i\la^{-1}s\Delta_g} \beta(P/\la) F(x,s)ds  \right)\\
     &= \scaledsoll{t} f_j(x)
\end{aligned}
\end{equation}
where for each $j$ 
\begin{equation}
    f_j(x) = \int_0^{j} e^{-i\la^{-1}s\Delta_g} \beta(P/\la) F(x,s)ds.
\end{equation}
It is important to note that each $f_j$ is a function which only depends on $x$. It also follows from the dual version of the homogeneous Strichartz estimate on logarithmic intervals \eqref{i9} that
\begin{equation}\label{fjbound}
    \lpnorm{f_j}{2}{M} \lesssim \la^{\frac{1}{2}}\mixlpnorm{F}{2}{q_e'}{M \times [0,T]} = \la^{\frac12}.
\end{equation}

Now, for each $t \in I_j$, define $H_j(\cdot,t) = \scaledsoll{t}f_j$. Using the fact that $\sigma(0)=1$, it is straightforward to check that for each $j$, $H_j = \sigma_{\la}H_j$.
Furthermore, by the support property of $\tilde \beta$ we have, $H_j = \tilde{\beta}(P/\la) \sigma_{\la}H_j$. Now, fix
\begin{equation}\label{eta}
\eta\in C^\infty_0((-3,3)) \quad
\text{with } \, \, \eta(t)=1, \, \, \, |t| \le 2.
\end{equation}
Define the function $\tilde{H}_j(x,t)= \eta(t-j) H_j(x,t)$ and the operator
\begin{equation}
    R_{\la,j}h(x,t) = \left( 1 - \eta(t-j)\right)\scaledsoll{t}h(x), \  \text{for} \  h \in L^{2}(M).
\end{equation}
We will need the following lemma, whose proof we postpone to the end of this section.
\begin{lemma}\label{3.1} For each $j$ we have
    \begin{equation} \label{time interval remainder estimate}
        \lpnorm{B \sigma_{\la}R_{\la,j}f_j}{\infty}{M \times I_j}  \lesssim_N \la^{-N}
\end{equation}
for any $N \geq 0$.
\end{lemma}
\noindent Note that for each $(x,t) \in M \times I_j$
\begin{equation}
    B \Sl F(x,t) = B\sigma_{\la}H_j(x,t) = B\sigma_{\la}\tilde{H}_{j}(x,t) + B \sigma_{\la}R_{\la,j} f_j(x,t)
\end{equation}
and thus for each $(x,t) \in A^{-}_{j}$
\begin{equation}\label{s20}
\begin{aligned}
    |B\sigma_{\la}\tilde{H}_j(x,t)| &\leq  |B\Sl F(x,t)| + |B \sigma_{\la}R_{\la,j}f_{j}(x,t)| \\
    &\lesssim_N \la^{\frac{d}4+ \frac12 + \varepsilon_0} + \la^{-N} \lesssim \la^{\frac{d}4+ \frac12 + \varepsilon_0}.
\end{aligned}
\end{equation}

Set $\tilde \sigma_\la=B\sigma_\la$ and $\theta_0 = \la^{-\e_1}$ for some fixed $0 < \e_1 < \frac{1}{2}$ to be determined later. Consider the the collection of microlocal pseudodifferential operators defined in \cite[section 2.3]{huang2024strichartz}, $\{A^{\theta_0}_\nu\}_{\nu \in \theta_0\Z{2d-1}}$.
The key  estimate that we require, which follows from Proposition 2.3 in \cite{huang2024strichartz} is the following.
\begin{proposition}\label{locprop} For $\theta_0$ defined above we have
    \begin{equation}\label{b1}
    \mixlpnorm{\tilde \sigma_\la \tilde H_j}{2}{q_e}{A^{-}_{j}}
    \lesssim  \othermixlitbiglpnorm{\tilde \sigma_\la A^{\theta_0}_\nu \tilde H_j}{t}{2}{q_e}{\nu}{x}{q_e}{M \times I_j}
    +\la^{\frac1{2}-} \lpnorm{\tilde H_j}{2}{M\times\mathbb{R}}.
    \end{equation}
\end{proposition}
\noindent
Here $\la^{\frac1{2}-}$ means $\la^{\frac1{2}-\e}$ for some unspecified $\e > 0$ which may depend on $\e_0,\,\e_1$.

\eqref{b1} follows from applying Proposition 2.3 in \cite{huang2024strichartz} to $\la^{-1/2}\tilde H_j$, and using that
\begin{equation} \label{necessary properties for prop 3.1}
    \lpnorm{\la^{-\frac12}\tilde{\sigma}_{\la}\tilde{H_{j}}}{\infty}{A^{-}_{j}} \lesssim \la^{\frac{d}4+\e_0}.
\end{equation}
Actually, if we let $n=d+1$ be the total dimension of $M\times \R$,  the $L^{\infty}$ bound $\la^{\frac{n-1}4+\e_0}$ in \eqref{necessary properties for prop 3.1} is exactly the ceiling bound used in the proof of Proposition 2.3 in \cite{huang2024strichartz}.

Recall that free Schr\"odinger solutions \(e^{it\Delta}\) in Euclidean space are frequency supported near the paraboloid
$
(\xi,-|\xi|^2)$.
The operators \(A^{\theta_0}_\nu\) may be viewed as analogues of wave packets in Euclidean space. They involve localization in frequency to cubes of size \(\lambda\theta_0\), supported near the paraboloid, together with spatial localization to tubular neighborhoods of geodesics. We also note that for $\theta_0 = \la^{-\e_1}$, the number of different \(A^{\theta_0}_\nu\) operators is $O(\la^{(2d-1)\e_1})$. See \cite{huang2024strichartz} for more details. 

The operators $\{\AvtO\}$ satisfy the following estimates due to the almost orthogonality of the individual operators
\begin{equation}\label{almost orthogonality}
    \mixlitbiglpnorm{\AvtO h}{p}{\nu}{p}{M}\ \lesssim \lpnorm{h}{p}{M}, \,\,\ 2 \leq p \leq \infty.  
\end{equation}
By duality this implies that for any sequence of functions $\{f_{\nu}\} \subset L^{p}(M)$ we have
\begin{equation}\label{duala}
    \lpnorm{\sum_{\nu}\big(\AvtO\big)^{*}f_{\nu}}{p}{M} \lesssim \mixlitbiglpnorm{f_{\nu}}{p}{\nu}{p}{M}, \ \ 1 \leq p \leq 2.
\end{equation}
We also require the following estimate from \cite{huang2024strichartz} for the commutator of $\AvtO$ and $\sigma_{\la}$.
\begin{lemma}\label{comprop} 
If $\delta>0$ in \eqref{22.6} is small enough
and
$\theta_0 =\la^{-\e_1}$ we have
for 
$B$ as in \eqref{22.9}
\begin{equation}\label{cc3}
\bigl\| \, 
B \sigma_\la A^{\theta_0}_\nu -B A^{\theta_0}_\nu
\sigma_\la \, \bigr\|_{L^2_{x,t}\to{L^{2}_tL^{q_e}_x}}
=O(\la^{2\e_1}).
\end{equation}
\end{lemma}

Now, by \eqref{cc3} and the fact that $\ell^{q_e}\subset \ell^2$, we have 
\begin{equation}\label{22.60}
    \begin{aligned}
        \bigl\| &
 \tilde \sigma_\la A^{\theta_0}_\nu  \tilde H_j
\bigr\|_{L^{2}_t\ell_\nu^{q_e}L^{q_e}_x}
\\
&\lesssim  \bigl\| 
B A^{\theta_0}_\nu
\sigma_\la   \tilde H_j
\bigr\|_{L^{2}_t\ell_\nu^{q_e}L^{q_e}_x} +\| \, 
(B \sigma_\la A^{\theta_0}_\nu -B A^{\theta_0}_\nu
\sigma_\la )\tilde H_j\,\|_{\ell_\nu^{2}L^{2}_tL^{q_e}_x} 
\\
&\lesssim 
\bigl\| 
B A^{\theta_0}_\nu
\sigma_\la  \tilde H_j
\bigr\|_{L^{2}_t\ell_\nu^{q_e}L^{q_e}_x}+\la^{2\e_1} \| \, 
\tilde H_j\,\|_{\ell_\nu^{2}L^2_{t,x}} .
    \end{aligned}
\end{equation}
Using that the number of choices of $\nu$ in the sum in \eqref{b1} is $O(\la^{(2d-1)\e_1})$ and $\tilde H_j$ is independent of $\nu$, we have 
\begin{equation} \label{Hj l^2L^2 bound}
    \| \, 
\tilde H_j\,\|_{\ell_\nu^{2}L^2_{t,x}}\lesssim\la^{(d-\frac12)\e_1}\| \, 
\tilde H_j\,\|_{L^2_{t,x}}.
\end{equation} 
Thus, if we choose $\e_1<\tfrac{1}{2d+3}$,  the second term on the right side of \eqref{22.60} is bounded by 
$\la^{\frac1{2}-}\|\tilde H_j\|_{L^2_{t,x}}$.

On the other hand, a similar argument to the proof of Lemma 3.1 yields for each $(x,t) \in M \times I_j$,
\begin{equation}\label{s20a}
   B A^{\theta_0}_\nu
\sigma_\la  \tilde H_j=B A^{\theta_0}_\nu
\sigma_\la   H_j+O(\la^{-N})= B A^{\theta_0}_\nu
   H_j+O(\la^{-N})
\end{equation}
Note that by \eqref{fjbound}
\begin{equation} \label{Hj L^2 bound}
   \la^{\frac1{2}-} \lpnorm{\tilde H_j}{2}{M\times\mathbb{R}} \sim \la^{\frac1{2}-}\lpnorm{f_j}{2}{M} \lesssim \la^{1-}.
\end{equation}
Thus, plugging \eqref{Hj l^2L^2 bound} into  \eqref{22.60} and using \eqref{s20a} and \eqref{Hj L^2 bound} we have 
\begin{align}\label{s15}
\| B S_\la F
\|_{L^{2}_tL^{q_e}_x(A^{-}_j)}
\lesssim
\| 
B \Atn S_\la F\|_{L^{2}_t\ell_\nu^{q_e}L^{q_e}_x(M \times I_j)}\,
+O(\la^{1-})
\end{align}
Summing over $A^{-}_j$ for $0\le j\le [T]$ and using \eqref{22.10}, we have 
\begin{align}\label{s15a}
\| B S_\la F
\|_{L^{2}_tL^{q_e}_x(A^{-})}
\lesssim
\| 
\Atn S_\la F\|_{L^{2}_t\ell_\nu^{q_e}L^{q_e}_x(M\times [0,T])}\,
+O(\la^{1-}).
\end{align}
Therefore, the proof of \eqref{small height sb} reduces to the following proposition, whose proof we postpone to Section 5.
\begin{proposition}
    \begin{equation} \label{microlocalized strichartz}
    \othermixlitbiglpnorm{A_{\nu}^{\theta_0}\Sl F}{t}{2}{q_e}{\nu}{x}{q_e}{M \times [0,T]} \lesssim \la \mixlpnorm{F}{2}{q_e'}{M \times [0,T]}.
    \end{equation}
\end{proposition}

We now turn to the proof of Lemma~\ref{3.1}. First, let $\{\la_{\ell}\}_{\ell \geq 0}$ denote the eigenvalues of $P$, with eigenfunctions $\{e_{\ell}\}_{\ell \geq 0}$ and 
\begin{equation}
    E_{\ell}f(x) = \int_{M}f(y)e_{\ell}(x)\overline{e_{\ell}(y)}dy
\end{equation}
Note that 
\begin{equation}
    \sigma_{\la}R_{\la,j}f(x,t) = \sum_{\ell \geq 0}\psi_{\la,j}( \la_{\ell}; t)E_{\ell}f(x)
\end{equation}
where 
\begin{equation}
    \psi_{\la,j}(\la_{\ell};t) = \frac{1}{(2\pi)^2}\iiint\widehat{\sigma}(\rho)\tilde{\beta}(\tau/\la)(1 - \eta(s-j))e^{i\phi(x,t,s,\rho, \tau, \la, \la_{\ell})}d\tau d\rho ds
\end{equation}
is an oscillatory integral with a phase function given by
\begin{equation}
    \phi(x,t,s,\rho, \tau, \la, \la_{\ell}) = (t-s)\tau +(\la^{\frac{1}{2}}\tau^{\frac12} - \la_{\ell})\rho{\color{red}-}s\la^{-1}\la_{\ell}^{2}.
\end{equation}
Note that 
\begin{equation}
    \big|\partial_{\tau}\phi\big| = \Big|t-s + \frac12\la^{\frac12}\tau^{-\frac12}\rho\Big| \gtrsim 1 + |t-s|
\end{equation}
if $\delta$ in \eqref{22.6} is chosen sufficiently small, since $|\tau| \sim \la$, $t \in I_j$ and $|s-j| \geq 2$ by the support properties of $(1 - \eta(s-j))$. Also, for $k \geq 2$ we have $\partial^{k}_{\tau}\phi = O(\la^{-k+1})$. Thus, by non-stationary phase
\begin{equation}
    \int \tilde{\beta}(|\tau|/\la)e^{i\phi(x,t,s,\rho, \tau, \la, \la_{\ell})}d\tau = O\big(\la(1 + \la + \la|t-s|)^{-N}\big) \ \ \ \text{for any} \ \ N \geq 0
\end{equation}
which implies
\begin{equation} \label{time interval remainder term kernel estimate}
    \psi_{\la,j}(\la_{\ell};t) = O\big( \la^{-N}\big) \ \ \text{for any} \ \ N \geq 0.
\end{equation}
Then by Sobolev estimates, \eqref{22.10} and \eqref{time interval remainder term kernel estimate} we have
\begin{equation}
    \lpnorm{B\sigma_{\la}R_{\la,j}f(\cdot, t)}{\infty}{M} \lesssim \la^{\frac{d}{2}}\lpnorm{\sigma_{\la}R_{\la,j}f(\cdot, t)}{2}{M} \lesssim_N \la^{-N}\lpnorm{f}{2}{M}
\end{equation}
for each $t \in I_j$. Thus, recalling \eqref{fjbound}, we have
\begin{equation}
    \lpnorm{B \sigma_{\la}R_{\la,j}f_j}{\infty}{M \times [j,j+1)}  \lesssim_N \la^{-N}.
\end{equation}
This completes the proof of Lemma~\ref{3.1}.

\section{Large Height Estimate}
Now we turn to the proof of \eqref{large height sb}. Let $\1_{A^+}(x,t)$ denote the indicator function of the set $A^+$ and 
choose $g(x,t)$ such that

\begin{equation} \label{uisng duality in large height argument}
    \mixlpnorm{g}{2}{q_e'}{A^{+}} = 1
    \quad \text{and } \, \,
    \mixlpnorm{B \Sl F}{2}{q_e}{A^+}
    =\iint  B \Sl F(x,t)
    \overline{\bigl(\1_{A^+}\cdot g\bigr)(x,t)} \, dx dt. 
\end{equation}
We have 
\begin{equation}
\begin{split}
    \iint  B \Sl F(x,t)
    \overline{\bigl(\1_{A^+}\cdot g\bigr)(x,t)} \, dx dt =   \iint  F(x,s)
    \overline{\Sl^*B^{*}\bigl(\1_{A^+}\cdot g\bigr)(x,s)} \, dx ds,
\end{split}
\end{equation}
where 
$$ \Sl^* G(x,s) = \int_0^T \scaledsoll{(s-t)} \1_A(t,s)\beta(P/\la) G(x,t)dt.
$$
Further, using \eqref{22.8} and \eqref{22.9} we have
\begin{equation} \label{microlocal decomposition in large height argument}
    \iint  F(x,s)
    \overline{\Sl^*B^{*}\bigl(\1_{A^+}\cdot  g\bigr)(x,s)} \, dx ds = \sum_{k = 1}^{N_0}\iint  F(x,s)
    \overline{B_{k, \la}\Sl^*B^{*}\bigl(\1_{A^+}\cdot g\bigr)(x,s)} \, dx ds.
\end{equation}
Recall that $N_0$ does not depend on $\la$ or $T$, only on the manifold $M$. For each $k$ we introduce a second height decomposition. Define the sets
\begin{equation}
  C^{-}_k=\{(x,s)\in M \times [0,T]: |B_{k,\la}S_\la^*B^{*} \big(\1_{A_+}\cdot g\big)|\lesssim \la^{\frac{d}{4}+\frac12+\e_0} \}
\end{equation}
and 
\begin{equation}
  C^{+}_k=\{(x,s)\in M \times [0,T] : |B_{k,\la}S_\la^*B^{*} \big(\1_{A_+}\cdot g\big)|\gtrsim \la^{\frac{d}{4}+\frac12+\e_0} \}
\end{equation}

\noindent Now, for each $k$, write 
\begin{equation} \label{decomposition into small and large heights in large height argument}
\begin{split}
    \iint  F(x,s)
    \overline{B_{k, \la}\Sl^*B^{*}\bigl(\1_{A^+}\cdot g\bigr)(x,s)} \, dx ds &=  \iint \1_{C^{-}_k}(x,s)  F(x,s)
    \overline{B_{k, \la}\Sl^*B^{*}\bigl(\1_{A_+}\cdot g\bigr)} \, dx ds \\
    & \quad +  \iint \1_{C^{+}_k}(x,s)  F(x,s)\overline{B_{k, \la}\Sl^*B^{*}\bigl(\1_{A_+}\cdot g\bigr)(x,s)} \, dx ds \\
    &= I + II
\end{split}
\end{equation}
By H\"older's inequality, we have 
\begin{equation} \label{initial estimate for I in large height argument}
\begin{split}
    I \leq \mixlpnorm{F}{2}{q_e'}{M \times [0,T]}\mixlpnorm{B_{k, \la}\Sl^*B^{*}\bigl(g \cdot \1_{A^{+}} \bigr)}{2}{q_e}{C^{-}_k}.
\end{split}
\end{equation}
It is straightforward to check that we can repeat the arguments from Section 3 to show that 
\begin{equation} \label{small height argument redux in large height argument}
\mixlpnorm{B_{k, \la}\Sl^*B^{*}\bigl(\1_{A^{+}}\cdot g\bigr)}{2}{q_e}{C^{-}_{k}} \lesssim \la\mixlpnorm{B^{*}\bigl(\1_{A^{+}} \cdot g\bigr)}{2}{q_e'}{M \times [0,T]} \lesssim \la,
\end{equation}
where the second inequality follows from the dual version of \eqref{22.10}  and the first part of \eqref{uisng duality in large height argument}. Plugging this back into \eqref{initial estimate for I in large height argument} yields 
\begin{equation} \label{final estimate for I in large height argument}
    I \lesssim \la.
\end{equation}

It remains to estimate the second term. We have
\begin{equation}
 II = \int_M \int_0^T\int_0^T \scaledsoll{(t-s)} \1_A(t,s)\beta(P/\la)B^{*}_{k,\la}(\1_{C^{+}_k} \cdot F )(x,s)\overline{B^{*}\big(\1_{A^+}\cdot g\big)(x,t)}dsdtdx
\end{equation}
If $|t-s|\le 1$, we can use \eqref{s1a} and repeat the arguments used in \eqref{strichartz estimate for unit time error term} to get the desired bound. If $|t-s|\ge 1$, we may use the following kernel bound
\begin{equation} \label{kernel bound in large height argument}
    \bigl|\scaledsoll{(t-s)}\beta^2(P/\la)(x,y)\bigr| \lesssim \la^{\frac{d}{2}}|t-s|^{-\frac{d}{2}}\exp(C_M|t-s|), \quad \text{if} \ |t-s| \leq T.
\end{equation}
This follows from Proposition 4.1 of \cite{blair2023strichartz}. The proof is based on lifting the relevant calculations to the universal cover of $M$ and applying the Hadamard parametrix on the covering manifold.
One may view this as a dispersive estimate with exponential loss up to the Ehrenfest time. By \eqref{kernel bound in large height argument}, we have
\begin{equation}
\begin{aligned}
\int_0^T&\int_0^T \int_{M} \scaledsoll{(t-s)} \1_A(t,s)\beta(P/\la)B^{*}_{k, \la}(\1_{C^{+}_k} \cdot F )(x,s)\overline{B^{*}\big(\1_{A^+}\cdot g\big)(x,t)}dxdsdt \\ 
&\lesssim \la^{\frac{d}2} e^{CT} \|B^{*}_{k, \la}\bigl(\1_{C^{+}_k} \cdot F\bigr) \|_{L^1_{x,s}}\|B^{*}\bigl(\1_{A^+} \cdot g\bigr) \|_{L^1_{x,t}} \\
& \lesssim \la^{\frac{d}2} e^{CT} \|F\|_{L^2_sL^{q_e'}_x(M \times [0,T])} 
\|\1_{C^{+}_k}\|_{L^2_sL^{q_e'}_x(M \times [0,T])} 
\mixlpnorm{g}{2}{q_e'}{M \times [0,T]}\mixlpnorm{\1_{A^+}}{2}{q_e}{M \times [0,T]} \\
& \lesssim \la^{\frac{d}2} e^{CT} \la^{-\frac{d}4-\frac12-\e_0}\norm{B_{k,\la}S_\la^*B^{*} \big(\1_{A^+}\cdot g\big)}_{L_{s}^{2}L_x^{q_e}(M \times [0,T])} \cdot \la^{-\frac{d}4-\frac12-\e_0}\mixlpnorm{B S_\la F}{2}{q_e}{A^{+}}.
\end{aligned}
\end{equation}
Here the second inequality follows from the dual version of \eqref{22.10} and H\"{o}lder's inequality. The last inequality follows from the definitions of the sets $A^{+}$ and $C^{+}_{k}$.

Fix \(I_j=(j-1,j]\) for \(1\le j\le T\). Then, for fixed \(s\in I_j\), the condition \(t\in \supp \1_A(t,s)\) is equivalent to \(t\in [j,T]\). Thus, by using the homogeneous Strichartz estimate \eqref{i9} and its dual version, we have
\begin{equation}\label{crude}
\begin{split}
    &\norm{B_{k,\la}S_\la^*B^{*} \big(\1_{A^+}\cdot g\big)}_{L_{s}^{2}L_x^{q_e}(M \times [0,T])} = \Big( \sum_{1\le j\le T}\norm{B_{k,\la}S_\la^*B^{*} \big(\1_{A^+}\cdot g\big)}_{L_{s}^{2}L_x^{q_e}(M \times I_j)} ^2\Big)^{\frac{1}{2}} \\
    &= \bigg( \sum_{1\le j\le T}\norm{B_{k,\la}\scaledsoll{s}\tilde{\beta}(P/\la)\int_{j}^{T}e^{-i\la^{-1}t\Delta_g}\beta(P/\la)B^{*}\bigl(\1_{A^+}\cdot g\bigr)(\cdot, t)dt}_{L_{s}^{2}L_x^{q_e}(M \times I_j)} ^2\bigg)^{\frac{1}{2}} \\
    &\lesssim \la^{\frac12}\bigg( \sum_{1 \leq j \le T}\lpnorm{\int_{j}^{T}e^{-i\la^{-1}t\Delta_g}\beta(P/\la) B^{*}\bigl(\1_{A^+}\cdot g\bigr)(\cdot, t)dt}{2}{M}^2\bigg)^{\frac{1}{2}} \\
    &\lesssim \la T^{\frac12} \mixlpnorm{B^{*}\bigl(\1_{A^{+}}\cdot g\bigr)}{2}{q_e'}{M \times [0,T]} \lesssim  \la T^{\frac12} \mixlpnorm{\1_{A^{+}}\cdot g}{2}{q_e'}{M \times [0,T]} \leq \la T^{\frac12}.
\end{split}
\end{equation}
By \eqref{crude}, if we choose $c_0$ sufficiently small (depending on $\e_0$) in the definition of $T$, so that $C e^{CT}T^{\frac12}\la^{-2\e_0} = C\la^{Cc_0}(c_0\log \la)^{\frac12}\la^{-2\e_0} \leq \frac1{2N_0}$, we have 
\begin{equation} \label{estimate for II in large height argument}
\begin{aligned}
II  \le \frac1{2N_0} \mixlpnorm{B \Sl F}{2}{q_e}{A^{+}}.
\end{aligned}
\end{equation}
Plugging \eqref{estimate for II in large height argument} and \eqref{final estimate for I in large height argument} back into \eqref{decomposition into small and large heights in large height argument} yields 
\begin{equation}
    \iint  F(x,s)
    \overline{B_{k, \la}\Sl^*B^{*}\bigl(\1_{A^+}\cdot g\bigr)(x,s)} \, dx ds \leq C\la + \frac{1}{2N_0}\mixlpnorm{B \Sl F}{2}{q_e}{A^{+}}.
\end{equation}
Finally, using this estimate and summing over $k$ in \eqref{microlocal decomposition in large height argument} yields
\begin{equation}
    \mixlpnorm{B \Sl F}{2}{q_e}{A^{+}} \leq C\la + \frac{1}{2}\mixlpnorm{B \Sl F}{2}{q_e}{A^{+}}
\end{equation}
This finishes the proof of \eqref{large height sb} and thus completes the proof of \eqref{s4}.

\section{Proof of Proposition 3.4}
Recall that we want to show 
\begin{equation} \label{microlocalized endpoint to endpoint estimate a}
    \othermixlitbiglpnorm{\AvtO\Sl F}{t}{2}{q_e}{\nu}{x}{q_e}{M \times [0,T]} \lesssim \la \mixlpnorm{F}{2}{q_e'}{M \times [0,T]}
\end{equation}
where
\begin{equation} \label{microlocalized inhomogeneous solution a}
    \AvtO\Sl F(x,t) = \int_{0}^{t}\AvtO \scaledsoll{(t-s)}{\1_A(t,s)}\beta(P/\la)F(\cdot, s)(x)ds.
\end{equation}
Note that by \eqref{almost orthogonality} and \eqref{strichartz estimate for unit time error term}, if we define
\begin{equation} \label{microlocalized inhomogeneous solution}
    \Sl^0 F(x,t) = \int_{0}^{t} \scaledsoll{(t-s)} \beta(P/\la)F(\cdot, s)(x)ds
\end{equation}
then 
\begin{equation} \label{microlocalized endpoint to endpoint estimate error}
    \othermixlitbiglpnorm{\AvtO(\Sl^0-\Sl) F}{t}{2}{q_e}{\nu}{x}{q_e}{M \times [0,T]} \lesssim  \mixlpnorm{(\Sl^0-\Sl) F}{2}{q_e}{M \times [0,T]}
    \lesssim \la \mixlpnorm{F}{2}{q_e'}{M \times [0,T]}
\end{equation}
Thus, to prove \eqref{microlocalized endpoint to endpoint estimate a}, it suffices to show that
\begin{equation} \label{microlocalized endpoint to endpoint estimate}
    \othermixlitbiglpnorm{\AvtO\Sl^0 F}{t}{2}{q_e}{\nu}{x}{q_e}{M \times [0,T]} \lesssim \la \mixlpnorm{F}{2}{q_e'}{M \times [0,T]}.
\end{equation}
We need the following dispersive type estimate
\begin{equation} \label{microlocalized dispersive estimate}
    \mixlitbiglpnorm{\AvtO \scaledsoll{(t-s)}\beta^2(P/\la)f}{\infty}{\nu}{\infty}{M} \lesssim \la^{\frac{d}{2}}|t-s|^{-\frac{d}{2}}\lpnorm{f}{1}{M}, \ \ |t-s| \lesssim T.
\end{equation}
\eqref{microlocalized dispersive estimate} follows from the proof of  \cite[Proposition 4.2]{blair2023strichartz}. Note that in \cite{blair2023strichartz}, the parameter is chosen as $\theta=\la^{-1/8}$. However, the same proof of the microlocalized kernel estimates extends to the larger scale $\theta_0=\la^{-\e_1}$ for arbitrarily small $\e_1>0$, provided that we take $T=c_0\log\la$ with some sufficiently small constant $c_0\ll \e_1$. 
The main idea here is to exploit the fact that the microlocalization in phase space provided by the $\AvtO$ operators restricts the number of ways a fixed point $x$ can propagate to $y$ along geodesics over a time interval of length $\approx |t-s|$.

\begin{remark}
    Note that if the operator was of the form $\AvtO e^{i\la^{-1}(t-s)\laplac{g}}\beta^2(P/\la)\big(\AvtO\big)^{*}$ in \eqref{microlocalized endpoint to endpoint estimate} and \eqref{microlocalized dispersive estimate} , then we would be able to directly apply the abstract Keel-Tao Theorem from \cite{KT} to prove that the estimate \eqref{microlocalized endpoint to endpoint estimate} holds. Nevertheless, we can still follow the proof of Keel-Tao to prove (\ref{microlocalized endpoint to endpoint estimate}).
\end{remark}

By a standard duality argument it suffices to show that if we consider a sequence of functions $\{\phi_{\nu}\}$ on $[0,T]\times M$ and define the bilinear operator 
\begin{equation}
    \bilinKT{}{F} = \int_{0}^{T}\int_{0}^{t}\sum_{\nu}\int_{M} \AvtO e^{i\la^{-1}(t-s)\laplac{g}}\beta(P/\la)(F(\cdot, s))(x)(x)\overline{\phi_{\nu}(x,t)}dxdsdt
\end{equation}
then we have
\begin{equation}
    \big| \bilinKT{}{F}\big| \lesssim \la \mixlpnorm{F}{2}{q_e'}{M \times [0,T]}\othermixlitbiglpnorm{\phi_{\nu}}{t}{2}{q_e'}{\nu}{x}{q_e'}{M \times [0,T]}.
\end{equation}
First, note that by a simple rescaling argument, \eqref{i9} is equivalent to
\begin{equation}
    \mixlpnorm{e^{i\la^{-1}t\laplac{g}}\beta(P/\la)f}{p}{q}{M \times [0,T]} \lesssim \la^{\frac{1}{p}}\lpnorm{f}{2}{M}
\end{equation}
whenever $(p,q)$ are admissible exponents, i.e. whenever we have
\begin{equation}
    d\Big(\frac{1}{2} - \frac{1}{q}\Big) = \frac{2}{p}, \ \ 2 \leq p \leq \infty, \ \ d \geq 3.
\end{equation}The dual version of this estimate says that
\begin{equation}\label{dual}
    \lpnorm{\int_{0}^{T}e^{-i\la^{-1}s\laplac{g}}\beta(P/\la)F(\cdot, s)ds}{2}{M} \lesssim \la^{\frac{1}{p}} \mixlpnorm{F}{p'}{q'}{M \times [0,T]}
\end{equation}
whenever $(p,q)$ are admissible.

To start we decompose our bilinear operator into dyadic intervals in the $s$ integral and consider the localized version
\begin{equation}
    \bilinKT{,j}{F} = \int_{0}^{T}\int_{t-2^{j+1}}^{t-2^j}\sum_{\nu}\int_{M} \AvtO e^{i\la^{-1}(t-s)\laplac{g}}\beta(P/\la)(F(\cdot, s))(x)\overline{\phi_{\nu}(x,t)}dxdsdt.
\end{equation}
Now, define 
\begin{equation}
    \gamma(a,b) = \frac{d}{2} - 1 - \frac{d}{2a} - \frac{d}{2b}.
\end{equation}
The following Lemma will be key to our argument. The reader should note that its role is the same as that of Lemma 4.1 in \cite{KT}.
\begin{lemma} \label{lemma 5.1}For any $\big(\frac{1}{a}, \frac{1}{b}\big)$ in a neighborhood of $\big( \frac{1}{q_e}, \frac{1}{q_e}\big)$ we have
    \begin{equation} \label{critical KT lemma}
        |\bilinKT{,j}{F}| \lesssim \la^{1 + \gamma(a, b)}2^{-\gamma(a,b)j}\mixlpnorm{F}{2}{a'}{[0,T] \times M} \othermixlitbiglpnorm{\phi_{\nu}}{t}{2}{b'}{\nu}{x}{b'}{[0,T] \times M}.
    \end{equation}
\end{lemma}
\noindent Lemma 5.1 follows from interpolating between the following special cases of (\ref{critical KT lemma}): \\
\begin{itemize}
    \item[(i)] $a = b =\infty$ \\
    \item[(ii)] $2 \leq a \leq q_e, \ \  b= 2$\\
    \item[(iii)] $a = 2, \ \ 2\leq b \leq q_e$.
\end{itemize}

We start by proving (i). For $k\in \mathbb{Z}$, let $I_k=[k2^j,(k+1)2^j)\cap [0,T]$,  $\tilde I_k=[(k-2)2^j, k2^j]\cap [0,T]$ and $I'_k = [(k+1)2^{j},(k+3)2^{j}) \cap [0,T]$. Now, note that 
\begin{equation}
    \gamma(\infty, \infty) = \frac{d}{2} - 1, \ \ \text{and} \ 1 + \gamma(\infty, \infty) = \frac{d}{2}.
\end{equation}
Also, note that for $s \in [t-2^{j+1}, t-2^j]$ we have $|t-s| \sim 2^j$. Thus, using (\ref{microlocalized dispersive estimate}) followed by Cauchy-Schwarz in both the $t$ and $s$ integrals before finally using Cauchy-Schwarz in the sum over $k$ and the finite overlap of the intervals $\{\tilde I_k\}$ we see that 
\begin{equation}
\begin{split}
    |\bilinKT{,j}{&F}| \\
    &\leq \sum_{k}\int_{I_k}\int_{t-2^{j+1}}^{t-2^{j}}\mixlitbiglpnorm{\AvtO e^{i\la^{-1}(t-s)\laplac{g}}\beta(P/\la)F(\cdot,s)}{\infty}{\nu}{\infty}{M}\sum_{\nu} \int_{M}|\phi_{\nu}(x,t)|dxdsdt \\
    &\lesssim \la^{\frac{d}{2}}2^{-\frac{d}{2}j}\sum_{k}\int_{I_k}\int_{t-2^{j+1}}^{t-2^{j}}\lpnorm{F(\cdot,s)}{1}{M}\mixlitbiglpnorm{\phi_{\nu}(\cdot, t)}{1}{\nu}{1}{M}dsdt \\
    &\lesssim \la^{\frac{d}{2}}2^{(1-\frac{d}{2})j}\sum_{k} \mixlpnorm{F}{2}{1}{M \times \tilde I_k}\othermixlitbiglpnorm{\phi_{\nu}}{t}{2}{1}{\nu}{x}{1}{M \times I_k} \\
    & \lesssim \la^{\frac{d}{2}}2^{(1-\frac{d}{2})j} \mixlpnorm{F}{2}{q_e}{M \times [0,T]} \othermixlitbiglpnorm{\phi_{\nu}}{t}{2}{1}{\nu}{x}{q_e}{M \times [0,T]}.
\end{split}
\end{equation}

Now we prove (ii). Fix $a \in [2, q_e]$. Then choose $p(a)$ so that $(p(a), a)$ is admissible. Note that 
\begin{equation}
    \gamma(a, 2) = -\frac{1}{p(a)'}, \ \ \text{and} \ 1 + \gamma(a,2) = \frac{1}{p(a)}.
\end{equation}
We start by noting that 
\begin{equation}
\begin{split}
    |\bilinKT{,j}{&F}| \\
    &\leq \sum_k \Big| \int_{I_k}\sum_{\nu}\int_{M} \AvtO e^{i\la^{-1}t\laplac{g}}\int_{t-2^{j+1}}^{t-2^j}e^{-i\la^{-1}s\laplac{g}}\beta(P/\la)(F(\cdot, s))(x)ds\overline{\phi_{\nu}(x,t)}dxdt\Big|.
\end{split}
\end{equation}
On each $I_k$ we apply  Cauchy-Schwarz to the integral over $M$ and the sum over $\{\nu\}$
\begin{equation}
    \begin{split}
    \Big| \int_{I_k}&\sum_{\nu}\int_{M} \AvtO e^{i\la^{-1}t\laplac{g}}\int_{t-2^{j+1}}^{t-2^j}e^{-i\la^{-1}s\laplac{g}}\beta(P/\la)(F(\cdot, s))(x)ds\overline{\phi_{\nu}(x,t)}dxdt\Big| \\
    &\leq \int_{I_k} \mixlitbiglpnorm{\AvtO e^{i\la^{-1}t\laplac{g}}\int_{t-2^{j+1}}^{t-2^j}e^{-i\la^{-1}s\laplac{g}}\beta(P/\la)(F(\cdot, s))(x)ds}{2}{\nu}{2}{M}\mixlitbiglpnorm{\phi_{\nu}(\cdot, t)}{2}{\nu}{2}{M}dt.
    \end{split}
\end{equation}
Then for any $t \in I_k$ we apply \eqref{almost orthogonality}, then use the fact that $\scaledsoll{t}$ is a unitary operator before finally using  \eqref{dual} to see that
\begin{equation}
\begin{split}
        \Bigg \lVert \AvtO \scaledsoll{t}\int_{t-2^{j+1}}^{t-2^j}&\scaledsol{s}\beta(P/\la)F(\cdot, s)ds \Bigg \rVert_{\ell_{\nu}^{2}L^{2}(M)} \\
        &\lesssim \lpnorm{ e^{i\la^{-1}t\laplac{g}}\int_{t-2^{j+1}}^{t-2^j}e^{-i\la^{-1}s\laplac{g}}\beta(P/\la)F(\cdot, s)ds}{2}{M} \\  
        &= \lpnorm{\int_{t-2^{j+1}}^{t-2^j}e^{-i\la^{-1}s\laplac{g}}\beta(P/\la)F(\cdot, s)ds}{2}{M} \\
        &\leq \sup_{t\in I_k}\lpnorm{\int_{t-2^{j+1}}^{t-2^j}e^{-i\la^{-1}s\laplac{g}}\beta(P/\la)F(\cdot, s)ds}{2}{M} \\
        &\lesssim \la^{\frac{1}{p(a)}} \mixlpnorm{F}{p(a)'}{a'}{M \times \tilde I_k}. 
\end{split}{}
\end{equation}
Thus, using the above estimate and then applying H\"{o}lder's inequality in each of the time integrals we have

\begin{equation}
    \begin{split}
    |\bilinKT{,j}{F}| & \lesssim \sum_{k}\left(\la^{\frac{1}{p(a)}} \mixlpnorm{F}{p(a)'}{a'}{M \times \tilde I_k} \right) \int_{I_k}\mixlitbiglpnorm{\phi_{\nu}(\cdot, t)}{2}{\nu}{2}{M}dt \\
    & \lesssim \sum_{k}\la^{\frac{1}{p(a)}} 2^{\frac{j}{p(a)'}}  \mixlpnorm{F}{2}{a'}{M \times \tilde{I_k}} \othermixlitbiglpnorm{\phi_{\nu}}{t}{2}{2}{\nu}{x}{2}{M \times I_k} \\
    &\lesssim \la^{\frac{1}{p(a)}} 2^{\frac{j}{p(a)'}}  \mixlpnorm{F}{2}{a'}{M \times [0,T]} \othermixlitbiglpnorm{\phi_{\nu}}{t}{2}{2}{\nu}{x}{2}{M \times [0,T]}.
    \end{split}
\end{equation}
Again, in the last inequality we used Cauchy-Schwarz in the sum over $k$ along with the finite overlap of the intervals $\{\tilde I_k\}$.

Finally, we prove (iii). Fix $b \in [2, q_e]$ and choose $p(b)$ so that $(p(b) , b)$ is admissible. Note that 
\begin{equation}
    \gamma(2, b) = -\frac{1}{p(b)'}, \ \ \text{and} \ 1 + \gamma(2,b) = \frac{1}{p(b)}.
\end{equation}
First, note that 
\begin{equation}
    \begin{split}
    |&\bilinKT{,j}{F}| \\
    &\leq \sum_k\Big| \int_{ I_k}\int_{M} F(y,s)\int_{[s+2^{j}, s+2^{j+1}]\cap [0,T]}\overline{e^{i\la^{-1}(s-t)\laplac{g}}\beta(P/\la)\sum_{\nu}\big(\AvtO \big)^{*}(\phi_{\nu}(\cdot,t))(y)}dtdyds\Big|.
    \end{split}
\end{equation}
On each ${I_k}$ we apply Cauchy-Schwarz to the integral over $M$ which yields
\begin{equation}
    \begin{split}
        \Big| &\int_{ I_k}\int_{M} F(y,s)\int_{[s+2^{j}, s+2^{j+1}]\cap [0,T]}\overline{e^{i\la^{-1}(s-t)\laplac{g}}\beta(P/\la)\sum_{\nu}\big(\AvtO \big)^{*}(\phi_{\nu}(\cdot,t))(y)}dtdyds\Big| \\
        &\leq \int_{ I_k} \lpnorm{F(\cdot, s)}{2}{M}\lpnorm{\int_{[s+2^{j}, s+2^{j+1}]\cap [0,T]}e^{i\la^{-1}(s-t)\laplac{g}}\beta(P/\la)\sum_{\nu}\big(\AvtO \big)^{*}(\phi_{\nu}(\cdot,t))dt}{2}{M}ds.
    \end{split}
\end{equation}
For any $s \in {I_k}$ we can first use the fact that $\scaledsol{s}$ is a unitary operator, then use \eqref{dual} followed by \eqref{duala}, which yields
\begin{equation}
    \begin{split}
\Bigg \lVert &\int_{[s+2^{j}, s+2^{j+1}]\cap [0,T]}\scaledsoll{(s-t)}\beta(P/\la)\sum_{\nu}\big(\AvtO \big)^{*}(\phi_{\nu}(\cdot,t))dt \Bigg \rVert_{L^{2}(M)} \\
&= \lpnorm{\scaledsoll{s}\int_{[s+2^{j}, s+2^{j+1}]\cap [0,T]}\scaledsol{t}\beta(P/\la)\sum_{\nu}\big(\AvtO \big)^{*}(\phi_{\nu}(\cdot,t))dt}{2}{M} \\
& = \lpnorm{\int_{[s+2^{j}, s+2^{j+1}]\cap [0,T]}\scaledsol{t}\beta(P/\la)\sum_{\nu}\big(\AvtO \big)^{*}(\phi_{\nu}(\cdot,t))dt}{2}{M} \\
&\leq \sup_{s \in  I_k}\lpnorm{\int_{[s+2^{j}, s+2^{j+1}]\cap [0,T]}\scaledsol{t}\beta(P/\la)\sum_{\nu}\big(\AvtO \big)^{*}(\phi_{\nu}(\cdot,t))dt}{2}{M} \\
&\lesssim \la^{\frac{1}{p(b)}}\mixlpnorm{\sum_{\nu}\big(\AvtO\big)^{*}\phi_{\nu}}{p(b)'}{b'}{M \times I'_k} \lesssim \la^{\frac{1}{p(b)}} \othermixlitbiglpnorm{\phi_{\nu}}{t}{p(b)'}{b'}{\nu}{x}{b'}{M \times I'_k}
    \end{split}
\end{equation}
Combining the above estimates and then applying H\"{o}lder's inequality in both of the time integrals yields
\begin{equation}
\begin{split}
    |\bilinKT{,j}{F}|
    &\lesssim \sum_k \left(\int_{I_k} \lpnorm{F(\cdot, s)}{2}{M}ds\right) \la^{\frac{1}{p(b)}} \othermixlitbiglpnorm{\phi_{\nu}}{t}{p(b)'}{b'}{\nu}{x}{b'}{M \times I'_k} \\
    &\lesssim\sum_k \la^{\frac{1}{p(b)}} 2^{\frac j {p(b)'}}\lpnorm{F}{2}{M \times I_k} \othermixlitbiglpnorm{\phi_{\nu}}{t}{2}{b'}{\nu}{x}{b'}{M \times I'_k} \\
    & \lesssim \la^{\frac{1}{p(b)}} 2^{\frac{j}{p(b)'}}\lpnorm{F}{2}{M \times [0,T]} \othermixlitbiglpnorm{\phi_{\nu}}{t}{2}{b'}{\nu}{x}{b'}{M \times [0,T]}.
\end{split}
\end{equation}
Once again, in the last inequality we used Cauchy-Schwarz in the sum over $k$ and the finite overlap of the intervals $\{I'_k\}$. This completes the proof of Lemma~\ref{lemma 5.1}.

Now, define a function on $M \times \theta_0\Z{2d-1} \times \R$ by $G(x,\nu,t) = \phi_{\nu}(x,t)$. Note that $M \times \theta_{0}\Z{2d - 1}$ is a measure space under the product measure of the volume form on $M$ and counting measure on $\theta_{0}\Z{2d-1}$. For each $t \in [0,T]$ the function $G(\cdot, t)$ lives on the measure space $M \times \theta_{0}\Z{2d-1}$ and the $L^{p}(M \times \theta_{0}\Z{2d-1})$ norm under this product measure is equivalent to the norm $\ell_{\nu}^{p}L^{p}(M)$. Thus, using Lemma 5.1 in \cite{KT} on the measure spaces $M$ and $M \times \theta_{0}\Z{2d-1}$ we can write
\begin{equation}\label{atomic decomp}
    F(x,t) = \sum_{\ell}f_{\ell}(t)\chi_{\ell}(x,t), \ \ \ G(x,\nu, t) = \sum_{k}g_k(t)\psi_{k}(x,\nu, t).
\end{equation}
For each $t, \ell$, the function $\chi_{\ell}(\cdot, t)$ is bounded by $O\big(2^{-\ell/q_e'}\big)$ and supported in a set with measure $O( 2^{\ell})$. Similarly, for each $t, k$, the function $\psi_{k}(\cdot, t)$ is bounded by $O\big(2^{-k/q_e'}\big)$ and supported in a set with measure $O(2^{k})$. Note that in this case, this is with respect to the aforementioned product measure on $M \times \theta_0\Z{2d-1}$. These decompositions also satisfy the crucial property that
\begin{equation} \label{crucial property of atomic decomposition}
\begin{split}
    \altlpnorm{\Big( \sum_{\ell}|f_{\ell}|^{q_e'}\Big)^{\frac{1}{q_e'}}}{2} &\lesssim \mixlpnorm{F}{2}{q_e'}{M \times [0,T]}, \\ 
    \text{and} \ \ &\altlpnorm{\Big( \sum_{k}|g_{k}|^{q_e'}\Big)^{\frac{1}{q_e'}}}{2} \lesssim \lVert G \rVert_{L^{2}_{t} L^{q_e'}_{x,\nu}(M \times \theta_0\Z{2d-1} \times [0,T])} = \othermixlitbiglpnorm{\phi_{\nu}}{t}{2}{q_e'}{\nu}{x}{q_e'}{M \times [0,T]}.
\end{split}
\end{equation}

Now, choose $ m \in \Z{}_{\geq 0}$ such that $2^m \sim \la$. Then (\ref{critical KT lemma}) can be restated as
\begin{equation}
    |\bilinKT{,j}{F}| \lesssim \la 2^{\gamma(a,b)(m - j)}\mixlpnorm{F}{2}{a'}{M \times [0,T]} \othermixlitbiglpnorm{\phi_{\nu}}{t}{2}{b'}{\nu}{x}{b'}{M \times [0,T]}.
\end{equation}
Combining this with \eqref{atomic decomp} and using the fact that 
\begin{equation}
    \gamma(a,b) = \frac{d}{2}\big( \frac{1}{q_e} - \frac{1}{a}\big) + \frac{d}{2}\big( \frac{1}{q_e} - \frac{1}{b}\big)
\end{equation}
yields
\begin{equation}
    |\bilinKT{,j}{F}| \lesssim \la \sum_{\ell}\sum_{k}2^{\big(\ell - \frac{d}{2}(j - m)\big)\big(\frac{1}{q_e} - \frac{1}{a}\big)}2^{\big(k - \frac{d}{2}(j - m)\big)\big(\frac{1}{q_e} - \frac{1}{b}\big)}\altlpnorm{f_{\ell}}{2}\altlpnorm{g_{k}}{2}.
\end{equation}
By Lemma 5.1 we know that this estimate holds for any $\big(\frac{1}{a}, \frac{1}{b}\big)$ in a neighborhood of $\big(\frac{1}{q_e}, \frac{1}{q_e}\big)$ so we can optimize the above estimate to see that for some fixed $\varepsilon > 0$ we have
\begin{equation}
    |\bilinKT{,j}{F}| \lesssim \la \sum_{\ell}\sum_{k}2^{-\varepsilon |\ell - \frac{d}{2}(j - m)|}2^{-\varepsilon |k - \frac{d}{2}(j - m)|}\altlpnorm{f_{\ell}}{2}\altlpnorm{g_{k}}{2}.
\end{equation}
Summing in $j$ yields
\begin{equation}
\begin{split}
    |\bilinKT{}{F}| \lesssim \la \sum_{\ell}\sum_{k}(1 + |\ell - k|)2^{-\varepsilon|\ell - k|}\altlpnorm{f_{\ell}}{2}\altlpnorm{g_{k}}{2}.
\end{split}
\end{equation}
Note that $(1+ |\cdot|)2^{-\varepsilon|\cdot|}$ is absolutely summable. Thus, by using the Cauchy-Schwarz inequality followed by Young's inequality and the fact that $\ell^{q_e'}\subset \ell^{2}$, before finally using \eqref{crucial property of atomic decomposition} we have
\begin{equation}
\begin{split}
    |\bilinKT{}{&F}| \lesssim \la \sum_{\ell}\sum_{k} (1 + |\ell - k|)2^{-\varepsilon|\ell - k|}\altlpnorm{f_{\ell}}{2}\altlpnorm{g_{k}}{2} \\
    &\leq \la \Big( \sum_{\ell}\altlpnorm{f_{\ell}}{2}^{2}\Big)^{\frac{1}{2}}\Big(\sum_{\ell}\Big(\sum_{k}(1+ |\ell - k|)2^{-\varepsilon|\ell - k|}\altlpnorm{g_k}{2}\Big)^2 \Big)^{\frac{1}{2}} \\
    &\lesssim \la \Big( \sum_{\ell}\altlpnorm{f_{\ell}}{2}^{2}\Big)^{\frac{1}{2}}\Big(\sum_{\ell}\altlpnorm{g_{k}}{2}^{2}\Big)^{\frac{1}{2}} = \la \altlpnorm{\Big(\sum_{\ell}|f_{\ell}|^2\Big)^{\frac{1}{2}}}{2}\altlpnorm{\Big(\sum_{k}|g_{k}|^2\Big)^{\frac{1}{2}}}{2} \\
    &\lesssim \la \altlpnorm{\Big(\sum_{\ell}|f_{\ell}|^{q_e'}\Big)^{\frac{1}{q_e'}}}{2}\altlpnorm{\Big(\sum_{k}|g_{k}|^{q_e'}\Big)^{\frac{1}{q_e'}}}{2} \lesssim \la \mixlpnorm{F}{2}{q_e'}{[0,T] \times M} \othermixlitbiglpnorm{\phi_{\nu}}{t}{2}{q_e'}{\nu}{x}{q_e'}{[0,T]\times M}.
\end{split}
\end{equation}
Thus, the proof of Proposition 3.4 is complete.

\section{Growth of Sobolev Norms for Cubic NLS on 3-dimensional Manifolds}

Before we describe how to adjust the argument in \cite{planchon2017growth} using the double-endpoint version of Theorem~\ref{thm} to prove Theorem~\ref{thm4}, we first show how the double-endpoint version of Theorem~\ref{thm4} and the endpoint version of \eqref{i9} imply an endpoint  Strichartz estimate over arbitrarily large time intervals with no derivative loss on the forcing term. To be more specific, assume 
$v$ is a solution to \begin{equation}\label{inho}
    i\partial_{t}v + \laplac{g}v = F, \ \ \  v(x,0) = v_0.
\end{equation}
Let $T > 0$. By using \eqref{i9} and \eqref{i10}, we have
\begin{equation}\label{n1}
\begin{aligned}
     \| &\beta(P/\la)v\|_{L^{2}_tL^{q_e}_x(M\times [0, T])}\\
     &\lesssim  \| \beta(P/\la)v(\cdot,0)\|_{L^{2}(M)}+ \| \beta(P/\la)v(\cdot,T)\|_{L^{2}(M)} \\
     &\quad+  (\la/\log \la)^{\frac12}\|\beta(P/\la)v\|_{L^{2}_{x,t}(M\times [0, T])}+ \|\beta(P/\la)F\|_{L^{2}_tL^{q_e'}_x(M\times [0, T])}.   
\end{aligned}
\end{equation}
\eqref{n1} follows from repeating the arguments used in \cite[Proposition 5.4]{bouclet2007strichartz}. In that paper, the authors divide the interval $[0,T]$ into subintervals of length $\la^{-1}$. The first two terms on the right-hand side come from the contributions from the first and last of these subintervals, while the Strichartz estimate of Burq, Gerard and Tzvetkov from \cite{bgtmanifold} is used on the other intervals before summing everything back up. In our case, we can divide the interval $[0,T]$ into subintervals of length $\la^{-1}\log \la$ and apply \eqref{i9} and \eqref{i10} when necessary.

Fix $\rho\in C_0^\infty$ such that $\rho=1$ on $[-1,1]$ on $\rho=0$ outside $[-2,2]$. By \eqref{n1} and Littlewood-Paley theory, for any $\la \gg 1$ we have the high frequency estimate
\begin{equation}\label{n2}
\begin{aligned}
     \| &(1-\rho(P/\la))v\|_{L^{2}_tL^{q_e}_x(M\times [0, T])}\\
     &\lesssim  \| v\|_{L_t^\infty L_x^{2}(M\times [0,T])}+  (\log \la)^{-\frac12}\|v\|_{L^{2}_tH^{1/2}_x(M\times [0, T])}+ \|F\|_{L^{2}_tL^{q_e'}_x(M\times [0, T])}.
\end{aligned}
\end{equation}
 For low frequencies, \eqref{n1} also implies
 \begin{equation}\label{n3}
\begin{aligned}
     \| &\rho(P/\la)v\|_{L^{2}_tL^{q_e}_x(M\times [0, T])}\\
     &\lesssim  \| v\|_{L_t^\infty L_x^{2}(M\times [0,T])}+  \la^{\frac12}\|v\|_{L^{2}_tL^{2}_x(M\times [0, T])}+ \|F\|_{L^{2}_tL^{q_e'}_x(M\times [0, T])}.   
\end{aligned}
\end{equation}
Combining \eqref{n2} and \eqref{n3} yields that for any $\la\gg 1$ we have
\begin{equation}\label{n4}
\begin{aligned}
     \| v\|_{L^{2}_tL^{q_e}_x(M\times [0, T])}
     &\lesssim  \| v\|_{L_t^\infty L_x^{2}(M\times [0,T])}+   (\log \la)^{-\frac12}\|v\|_{L^{2}_tH^{1/2}_x(M\times [0, T])} \\
     &+ \la^{\frac12}\|v\|_{L^{2}_tL^{2}_x(M\times [0, T])}+ \|F\|_{L^{2}_tL^{q_e'}_x(M\times [0, T])}.
\end{aligned}
\end{equation}
This is the Strichartz estimate we will use in the proof of Theorem~\ref{thm4}. Note that it is the analog of Proposition 2.2 in \cite{planchon2017growth}.
It is worth pointing out that although Proposition 2.2 in \cite{planchon2017growth} is stated for $T \in (0,1)$ the argument from \cite[Proposition 5.4]{bouclet2007strichartz} shows that the method continues to work for all $T > 0$. It is also worth noting that \cite[Proposition 2.2]{planchon2017growth} includes a derivative loss of $\varepsilon > 0$ on the analog of the first term in \eqref{n4}. However, this can be done away with by applying Littlewood-Paley theory to the estimate \eqref{n1} and then taking the sup norm over $[0,T]$ as opposed to taking the sup norm over $[0,T]$ in the analogous estimate to \eqref{n1} in \cite[Proposition 5.4]{bouclet2007strichartz}.

Now, fix $d = 3$ and note that $q_e = 6$. Let $u$ be a solution to \eqref{cubic NLS}. For now, we consider the case when $m = 2k$ for $k\in \N$. Also, let $\tau > 0$. Note that throughout the remainder of this section, the implicit constants are allowed to depend on the $H^{s}$ norms of the initial data for $s \leq m$. By repeating the argument used in the proof of \cite[Proposition 3.5]{planchon2017growth}, using \eqref{n4} in place of \cite[Proposition 2.2]{planchon2017growth}, and directly applying \cite[Proposition 3.3]{planchon2017growth} along with the estimate
\begin{equation}
    \lVert u\rVert_{L^{\infty}_{t}H^{s}_{x}(M \times [0,\tau])} \lesssim \lVert u\rVert_{L^{\infty}_{t}H^{2k}_{x}(M \times [0,\tau])}^{\frac{s-1}{2k-1}} \ \ \ \text{for} \ s \in [1,2k],
\end{equation}
which is (24) in \cite{planchon2017growth}, 
we get the following analog of (32) and (33) in \cite{planchon2017growth}:
\begin{equation}
   \begin{aligned}\label{32}
\|\partial_t^j u\|_{L^{2}_tL^{6}_x(M\times [0, \tau])}
&\lesssim \|u\|_{L^\infty_t H^{2k}(M\times [0, \tau])}^{\frac{2j-1}{2k-1}} + \sqrt{\tau} (\log\la)^{-\frac12}\|u\|_{L^\infty_t H^{2k}(M\times [0, \tau])}^{\frac{4j-1}{4k-2}}   \\&\quad +\sqrt{\tau} (\la^{\frac12}+1)\|u\|_{L^\infty_t H^{2k}(M\times [0, \tau])}^{\frac{2j-1}{2k-1}}, 
\end{aligned} 
\end{equation}

\begin{equation}\label{33}
  \begin{aligned}
\|\partial_t^j u\|_{L^{2}_tW^{1,6}_x(M\times [0, \tau])}
&\lesssim \|u\|_{L^\infty_t H^{2k}(M\times [0, \tau])}^{\frac{2j}{2k-1}} + \sqrt{\tau} (\log\la)^{-\frac12}\|u\|_{L^\infty_t H^{2k}(M\times [0, \tau])}^{\frac{4j+1}{4k-2}}   \\&\quad +\sqrt{\tau} (\la^{\frac12}+1)\|u\|_{L^\infty_t H^{2k}(M\times [0, \tau])}^{\frac{2j}{2k-1}}.
\end{aligned}  
\end{equation}
Then, using \eqref{32}, \eqref{33} and the results from Section 3 in \cite{planchon2017growth},
we can repeat the arguments in Section 5 of \cite{planchon2017growth} to show that for any $\tau > 0$ we have
\begin{equation}\label{mainapplication}
\begin{split}
    \|u(\cdot, \tau)\|^2_{H^m(M)} - \|u(\cdot,0)\|^2_{H^m(M)} 
    \lesssim & \ \tau (\log\la)^{-1} \|u\|^2_{L_{t}^\infty H_{x}^m(M \times [0,\tau])} \\
    &+ (1+\tau\la)\|u\|^\gamma_{L_{t}^\infty H_{x}^m(M \times [0,\tau])}
\end{split}
\end{equation}
where $\gamma\in (0,2)$ is some fixed constant. Note that this argument critically relies on the use of the ``modified energies", $\mathcal{E}_{2k}(u)$, which are defined in Section 3 of \cite{planchon2017growth}.

Now, fix $T > 0$ and let $\tau = \sqrt{T}$ and $\la=\exp(C_0\sqrt{T})$ for some large constant $C_0$. Then \eqref{mainapplication} becomes
\begin{equation}
\begin{split}
    \|u(\cdot, \tau)\|^2_{H^m(M)} - \|u(\cdot,0)\|^2_{H^m(M)} 
    \leq & \ \frac12 \|u\|^2_{L_{t}^\infty H_{x}^m(M \times [0,\tau])} \\
    &+ C\sqrt{T}\exp(C_0 \sqrt{T})\|u\|^\gamma_{L_{t}^\infty H_{x}^m(M \times [0,\tau])}.
\end{split}
\end{equation}
Since $\gamma<2$ this further implies that we have 
\begin{equation}\label{mainapplication1}
    \|u(\cd,\tau)\|^2_{H^m(M)} 
    \leq 2\|u(\cd,0)\|^2_{H^m(M)} +C\exp(C\sqrt{T})
\end{equation}
if $C \gg C_0$ is sufficiently large.
This then yields that for any $t > 0$ we have 
\begin{equation}
    \|u(\cd, t+\tau)\|^2_{H^m(M)} 
    \leq 2\|u(\cd,t)\|^2_{H^m(M)} +C \exp(C \sqrt{T}).
\end{equation}
If we iterate \eqref{mainapplication1} $\sim \sqrt{T}$ many times, we get 

\[
\sup_{0\leq t \leq T} \|u(\cd,t)\|_{H^m(M)} \leq C\exp(C\sqrt{T}),
\]
where $C = C(m, \|f\|_{H^m}) > 0$.

Now, we turn to the case where $m  = 2k+1$ for some $k \in \N$. Note that if we assumed $f \in H^{2k+2}(M)$ then the result would follow from the fact that the desired growth estimate holds for $H^{2k+2}$. The main point here is that we only assume that $f \in H^{2k+1}$. We can prove the result in this case by following the argument outlined above for the case of even $m$, replacing the use of the modified energies $\mathcal{E}_{2k}(u)$ by the use of the modified energies $\mathcal{E}_{2k+1}(u)$ defined in Section 6 of \cite{planchon2017growth}. We leave the details to the reader.

\section{Strichartz Estimates for Schr\"{o}dinger Operators with critically singular potentials}
In this section, we shall prove Theorem~\ref{thm2}. Our argument is similar to the one in \cite{HSSchro}. If $\beta$ is as in \eqref{i6}, let us define
``wider cutoffs" that we shall also use as follows
\begin{equation}\label{m.2}
\tb(s)=\sum_{|j|<10}\beta(2^{-j}s)\in C^\infty_0(2^{-10},
2^{10}).
\end{equation}
For future use, note that
\begin{equation}\label{m.3}
\tb(s)=1 \quad \text{on } \, \, (1/8,8).
\end{equation}
We also remind the reader that $H_V = -\laplac{g} + V(x)$ where $V \in L^{d/2}(M)$, $P_V = \sqrt{H_V}$ and $q_e = \frac{2d}{d-2}$. Set $\e=\e(\la) = (\log\la)^{-1}$. Note that to prove Theorem 1.2 it suffices to prove the endpoint estimate
\begin{equation}\label{i11 endpoint}
\| e^{-itH_V} \beta(P_V/\la) f\|_{L^2_tL^{q_e}_x(M
\times [0, \, (\la\e)^{-1}])}\lesssim \, \|f\|_{L^2(M)}, \quad \la\gg 1.
\end{equation}

Let $w$ be the solution to the scaled inhomogeneous Schr\"{o}dinger equation
\begin{equation}\label{w1}
    \big(i \la \partial_t + \laplac{g} + i \e \la\big)w(x,t) = F(x,t) \ , \quad w(x,0) = 0
\end{equation}
where we also assume that
\begin{equation}\label{w2}
    w(\cdot,t), \,\,F(\cdot,t)=0 \,\,\, \text{if} \,\,\,t\notin  [0,\e^{-1}].
\end{equation}
By Duhamel's formula we have
\begin{equation}
    w(x,t) = \frac{1}{i \la} \int_{0}^{t}\scaledsoll{(t-s)}e^{-\e(t-s)}F(x,s)ds.
\end{equation}
For later use, if $\1_+(s)=\1_{[0,+\infty)}(s)$ denotes the Heaviside function, we recall that
\begin{equation*}
(2\pi)^{-1}\int_{-\infty}^\infty \frac{e^{it\tau}}{i\tau +\e} \, d\tau = \1_+(t)e^{-\e t}.
\end{equation*}
Thus,
\begin{equation} \label{frequency localized w}
\begin{aligned}
    \tb(P/\la) w(x,t) &= \frac{1}{i \la} \int_{0}^{t}\scaledsoll{(t-s)}e^{-\e(t-s)}\tb(P/\la)F(x,s)ds \\
    &=(2\pi)^{-1}\int_0^{\e^{-1}}\int_{-\infty}^\infty
\frac{e^{i(t-s)\tau}}{-\la\tau+\Delta_g+i\e\la}
\tb(P/\la)F(x,s)\, d\tau ds.
\end{aligned}
\end{equation}

A simple rescaling argument shows that the homogeneous  Strichartz estimate \eqref{i9} from \cite{huang2024strichartz} is equivalent to the estimate 
\begin{equation} \label{rescaled i9}
\mixlpnorm{\scaledsoll{t}\beta(P/\la)f}{p}{q}{M \times[0,\log\la]} \lesssim \la^{\frac{1}{p}}\lpnorm{f}{2}{M}.
\end{equation}
Combining this with the inhomogeneous Strichartz estimate \eqref{i10} implies the following two estimates for solutions of \eqref{w1}.
\begin{proposition}\label{6.1} Suppose $w$ and $F$ satisfy \eqref{w1} and \eqref{w2}, and $\e=(\log\la)^{-1}$, we have 
    \begin{equation}
        \mixlpnorm{\tb(P/\la)w}{2}{q_e}{M \times [0,\e^{-1}]} \lesssim \mixlpnorm{F}{2}{q_e'}{M \times [0,\e^{-1}]}
    \end{equation}
    and 
    \begin{equation}
        \mixlpnorm{\tb(P/\la)w}{2}{q_e}{M \times [0,\e^{-1}]} \lesssim \la^{-\frac{1}{2}}\e^{-\frac{1}{2}} \lpnorm{F}{2}{M \times [0,\e^{-1}]}
    \end{equation}
\end{proposition}
\noindent The first of these two estimates follows directly from the double-endpoint inhomogeneous Strichartz estimate \eqref{a} and (\ref{frequency localized w}). The second estimate follows from  \eqref{rescaled i9} and (\ref{frequency localized w}). We have
\begin{equation*}
\begin{split}
    \Big\lVert \tb(P/\la)w &\Big\rVert_{L^{2}_{t}L^{q_e}_{x}(M \times [0,\e^{-1}])} = \la^{-1}\mixlpnorm{\int_{0}^{t}\scaledsoll{(t-s)}e^{-\e(t-s)}\tb(P/\la)F(\cdot,s)ds}{2}{q_e}{M \times [0,\e^{-1}]} \\
    &= \la^{-1}\mixlpnorm{\int_{0}^{\e^{-1}}e^{-\e(t-s)}\scaledsoll{(t-s)}\tb(P/\la)F(\cdot,s)\mathbbm{1}_{s\leq t}(s,t)ds}{2}{q_e}{M \times [0,\e^{-1}]} \\
    & \leq \la^{-1} \int_{0}^{\e^{-1}}\mixlpnorm{\scaledsoll{(t-s)}\tb(P/\la)F(\cdot,s)}{2}{q_e}{M \times [0,\e^{-1}]}ds \\
    &\lesssim \la^{-\frac{1}{2}}\int_{0}^{\e^{-1}}\lpnorm{e^{-i\la^{-1}s\laplac{g}}\tb(P/\la)F(\cdot,s)}{2}{M}ds \\
    &\lesssim \la^{-\frac{1}{2}} \int_{0}^{\e^{-1}} \lpnorm{F(\cdot,s)}{2}{M}ds \leq\la^{-\frac{1}{2}}\e^{-\frac{1}{2}}\lpnorm{F}{2}{M \times [0,\e^{-1}]}
\end{split}
\end{equation*}
where we used Minkowski's inequality to go from the second to the third line, \eqref{rescaled i9} to go from the third to the fourth line, before finally using the fact that $\scaledsol{s}$ is a unitary operator and Holder's inequality in the last line.

Proposition~\ref{6.1} further implies the following two estimates involving Littlewood-Paley operators in the time variable.
\begin{proposition}\label{6.2} 
Let  $\beta$ be as in \eqref{i6}.
Suppose $w$ and $F$ satisfy \eqref{w1} and \eqref{w2}, and $\e=(\log\la)^{-1}$. Then we have 
\begin{equation} \label{Strichartz implies time localized 1}
        \mixlpnorm{\beta(-D_t/\la)w}{2}{q_e}{M \times [0,\e^{-1}]} \lesssim \mixlpnorm{F}{2}{q_e'}{M \times [0,\e^{-1}]}
    \end{equation}
and
\begin{equation} \label{Strichartz implies time localized 2}
    \mixlpnorm{\beta(-D_t/\la)w}{2}{q_e}{M \times [0,\e^{-1}]} \lesssim \la^{-\frac{1}{2}}\e^{-\frac{1}{2}} \lpnorm{F}{2}{M \times [0,\e^{-1}]}.
\end{equation}
\end{proposition}
\noindent To prove this, 
we first note that the kernel of $\beta(-D_t/\la)$
is $O(\la (1+\la|t-t'|)^{-2})$. 
Therefore, one can use Young's inequality in the time variable to see that 
$$\|\beta(-D_t/\la) \tb(P/\la)w\|_{L^2_t L^{q_e}_x(M \times \R)}
\lesssim \| \tb(P/\la)w\|_{L^2_tL^{q_e}_x(M \times \R)}.
$$
Thus, by Proposition \ref{6.1},
it suffices to show that $\beta(-D_t/\la)(I-\tb(P/\la))w$ enjoys the bounds in 
\eqref{Strichartz implies time localized 1} and \eqref{Strichartz implies time localized 2}.

Recalling \eqref{frequency localized w}, this means that it suffices
to show that
\begin{multline}\label{m.32}
\Bigl\|
\int_0^1\int_{-\infty}^\infty
\frac{e^{i(t-s)\tau}}{-\la\tau-P^2+i\e\la}
\beta(-\tau/\la) \, 
\bigl(I-\tb(P/\la)\bigr)\, F(\cd,s)\, d\tau ds
\, \Bigr\|_{L^2_tL^{q_e}_x(M \times \R)}
\\
\lesssim
\|F\|_{L^{2}_tL^{q_e'}_x(M \times [0,\e^{-1}])}.
\end{multline}
as well as
\begin{multline}\label{m.33}
\Bigl\|
\int_0^1\int_{-\infty}^\infty
\frac{e^{i(t-s)\tau}}{-\la\tau-P^2+i\e\la}
\beta(-\tau/\la) \, 
\bigl(I-\tb(P/\la)\, F(\cd,s)\, d\tau ds
\, \Bigr\|_{L^2_tL^{q_e}_x(M \times \R)}
\\
\lesssim \la^{-\frac12}\e^{-\frac12}\|F\|_{L^2_{x,t}(M \times [0,\e^{-1}])},
\end{multline}
If we set
$$\alpha(t,s;\mu)
= \int_{-\infty}^\infty
\frac{e^{i(t-s)\tau}}{-\la\tau-\mu^2+i\e\la}
\beta(-\tau/\la) \, 
\bigl(1-\tb(\mu/\la)\bigr) \, d\tau,$$
then by \eqref{m.3} and the support
properties of $\beta$ we have
for $j=0,1,2$ that
$$
\la \,
\bigl|
\la^{j} \,
\partial_\tau^j \bigl((1-\tb(\mu/\la)\bigr)  \beta(-\tau/\la)
(-\la\tau-\mu^2+i\e\la)^{-1}\bigr)\, \bigr|
\lesssim \la (\mu^2+\la^2)^{-1}.$$
 This then implies that
\begin{equation}\label{abound}
  |\alpha(t,s;\mu)|\lesssim 
\la(1+\la|t-s|)^{-2} \cdot (\mu^2+\la^2)^{-1}  
\end{equation}
using a simple integration by parts argument.
If we use \eqref{abound}, then \eqref{m.32} and \eqref{m.33} follow from simple Sobolev estimates together with an application of Young's inequality in the time variable; see, for example, Lemmas 2.2 and 2.3 in \cite{HSSchro} for more details.

Let us now show how we can use Proposition~\ref{6.2} to prove Theorem~\ref{thm2}. First, we note that just as in the $V = 0$ case, the eigenvalues of $P_V$ are nonnegative, discrete, and tend to infinity. Let $\{\mu_j\}$ denote these eigenvalues with associated orthonormal basis of $L^2$-normalized eigenfunctions $\{e^V_j\}$. We let $E^{V}_j$ denote the projection onto the $j$-th eigenspace, i.e.
$$ E^V_jf(x) = \int_{M}f(y)e^{V}_{j}(x)\overline{e^{V}_{j}(y)}dy.$$
To prove \eqref{i11 endpoint}, it clearly suffices to show that
for $\lambda \gg 1$ we have the uniform bounds
\begin{multline}\label{d.1}
\Bigl\| \, \eta(\la \e t) e^{-itH_V}f_\la
\Bigr\|_{L^2_tL^{q_e}_x(M \times [0,\e^{-1}])}\le C\|f_\la\|_{L^2(M)},
\\
\text{if } \, \, \mathrm{spec } \, f_\la
\subset [9\la/10, 11\la/10] \quad
\text{and } \, \, \,
\eta\in C^\infty_0((0,1)) \quad
\text{is fixed}.
\end{multline}
Note that when we write that $\mathrm{spec }f_{\la} \subset [9\la/10, 11\la/10]$ we mean that $E^{V}_{j} f = 0$ if $\mu_j \notin [9\la/10,11\la/10]$. This is equivalent to showing that 
\begin{equation} \label{rescaled d.1}
    \mixlpnorm{w}{2}{q_e}{M \times [0,\e^{-1}]} \lesssim \la^{\frac{1}{2}}\altlpnorm{f_{\la}}{2}, \ \text{where} \ w(x,t) = \eta(\e t)e^{-i\la^{-1}tH_V}f_{\la}(x).
\end{equation}

Write 
\begin{equation}
    V = V_{\leq \ell} + V_{> \ell}
\end{equation}
where 
\begin{equation}
    V_{> \ell}(x) =\begin{cases} 
    &V(x), \ \text{if } |V(x)| > \ell \\    
    & \ \ 0  \ \ \ , \ \text{otherwise. }
    \end{cases} 
\end{equation}
Since $V \in L^{d/2}$ we have
\begin{equation}
    \lpnorm{V_{> \ell}}{d/2}{M} = \delta(\ell) \to 0 \ \text{as} \ \ell \to \infty.
\end{equation}
Also, from the definition of $V_{\leq \ell}$, it follows that
\begin{equation} \label{potential height bound}
    \lpnorm{V_{\leq \ell}}{\infty}{M} \leq \ell.
\end{equation}

Now, since $-H_V = \laplac{g} - V$ we can write

\begin{equation}
    \big(i \la \partial_t + \laplac{g} + i \e\la\big)w = \big(i \la \partial_t - H_V + i \e\la\big)w  + Vw 
    = \big(i \la \partial_t - H_V + i \e\la\big)w + V_{\leq \ell}w + V_{>\ell}w
\end{equation}
and $w(\cd,0) = 0$. Thus, we can further write
\begin{equation}
    w = \widetilde{w} + w_{\leq \ell} + w_{> \ell}
\end{equation}
where
\begin{equation}
    \big(i \la \partial_t + \laplac{g} + i \e\la\big)\widetilde{w} = \big(i \la \partial_t - H_V + i \e\la\big)w = \widetilde{F}, \ \ \widetilde{w}(\cd,0) = 0,
\end{equation}
\begin{equation}
    \big(i \la \partial_t + \laplac{g} + i \e\la\big)w_{\leq \ell} = V_{\leq \ell}w = F_{\leq \ell}, \ \ w_{\leq \ell}(\cd,0) = 0,
\end{equation}
and
\begin{equation}
    \big(i \la \partial_t + \laplac{g} + i \e\la\big)w_{> \ell} = V_{> \ell}w = F_{> \ell}, \ \ w_{> \ell}(\cd,0) = 0.
\end{equation}
Furthermore, note that since $w(x,t) = 0$ if $t \notin [0,\e^{-1}]$ the forcing terms $\widetilde{F}, F_{\leq \ell}$ and $F_{> \ell}$ also vanish if $t \notin [0,\e^{-1}]$. Thus, we have
\begin{equation}
\begin{split}
    \mixlpnorm{w}{2}{q_e}{M \times [0,\e^{-1}]} &\leq \mixlpnorm{\beta(-D_t/ \la)w}{2}{q_e}{M \times [0,\e^{-1}]} + \mixlpnorm{\big(I - \beta(-D_t/ \la)\big)w}{2}{q_e}{M \times [0,\e^{-1}]} \\
    &\leq  \mixlpnorm{\beta(-D_t/ \la)\widetilde{w}}{2}{q_e}{M \times [0,\e^{-1}]} + \mixlpnorm{\beta(-D_t/ \la)w_{\leq \ell}}{2}{q_e}{M \times [0,\e^{-1}]} \\
    & \ \ \ \ \ \ + \mixlpnorm{\beta(-D_t/ \la)w_{>\ell}}{2}{q_e}{M \times [0,\e^{-1}]} + C_N\la^{-N}\altlpnorm{f_{\la}}{2}.
\end{split}
\end{equation}
Here we used the fact that 
\begin{equation}\label{d.4}
\begin{split}
 \mixlpnorm{\big(I - \beta(-D_t/ \la)\big)w}{2}{q_e}{[0,\e^{-1}]\times M} 
    \leq  C_N\la^{-N}\altlpnorm{f_{\la}}{2}.
\end{split}
\end{equation} 
To see this, note that the Fourier transform
of $t\to \eta(t)e^{-it\la^{-1}\mu^2_j}$ is
$\Hat \eta(\tau+\mu_j^2/\la)$  and so
\begin{equation}
\bigl(I-\beta(-D_t/\la)\bigr)w(x,t)
= \sum_{\mu_j\in [9\la/10, 11\la/10]} a(t;\mu_j) E^V_jf_\la(x),
\end{equation}
where
\begin{equation}\label{d.5}
a(t;\mu)=(2\pi)^{-1}
\int_{-\infty}^\infty e^{it\tau}
\hat{\eta}(\tau+\mu^2/\la) \, 
\big(1-\beta(-\tau/\la)\bigr) \, d\tau.
\end{equation}
Note that by the support assumption of $\beta$ in \eqref{i6}, it is not hard to check
$$ |\hat{\eta} (\tau+\mu^2/\la) \, 
\big(1-\beta(-\tau/\la)\bigr)|\lesssim_{N_1,N_2} \la^{-N_1} (1+|\tau|)^{-N_2}\,\,\text{for any} \, N_1,N_2 \geq 0 \,\, \text{if}\,\,\mu\in [9\la/10, 11\la/10].
$$
Additionally, by \cite{BHSS} we have 
\begin{equation}\label{sobolev}
    \|u\|_{L^{q_e}(M)}\lesssim \|(I+H_V)^{1/2}u\|_{L^2(M)}.
\end{equation}
Thus, \eqref{d.4} follows from the spectral theorem together with the Sobolev estimate \eqref{sobolev}; see also Lemma 3.1 in \cite{HSSchro} for more details.

As a result, we see that \eqref{rescaled d.1} would follow from the following three inequalities:
\begin{equation}
    \mixlpnorm{\beta(-D_t/ \la)\widetilde{w}}{2}{q_e}{[0,\e^{-1}]\times M} \lesssim \la^{\frac{1}{2}}\altlpnorm{f_{\la}}{2},
\end{equation}
\begin{equation}
    \mixlpnorm{\beta(-D_t/ \la)w_{\leq \ell}}{2}{q_e}{M \times [0,\e^{-1}]} \lesssim \ell \la^{\frac{1}{2}}\altlpnorm{f_{\la}}{2}
\end{equation}
and 
\begin{equation} \label{time localized endpoint to endpoint estimate for large height potential}
    \mixlpnorm{\beta(-D_t/ \la)w_{>\ell}}{2}{q_e}{M \times [0,\e^{-1}]} \leq \frac{1}{2}\mixlpnorm{w}{2}{q_e}{M \times [0,\e^{-1}]}.
\end{equation}
Here $\ell$ is fixed and will be determined at the conclusion of the argument. 

Note that 
\begin{equation}
    \widetilde{F} = \big(i \la \partial_t - H_V + i \e\la\big)w = i\e\la\big( \eta'(\e t) + \eta(\e t)\big)e^{-i\la^{-1}tH_V}f_{\la}
\end{equation}
so 
\begin{equation}
    \widetilde{F} = i\e\la\big( \eta'(\e t) + \eta(\e t)\big)e^{-i\la^{-1}tH_V}f_{\la}
\end{equation}
We have that (\ref{Strichartz implies time localized 2}) implies
\begin{equation}
\begin{split}
    \Big\lVert \beta(-D_t/ \la)\widetilde{w} \Big\rVert_{L^{2}_{t}L^{q_e}_{x}(M \times [0,\e^{-1}])} \lesssim \la^{-\frac{1}{2}}\e^{-\frac{1}{2}}\altlpnorm{\widetilde{F}}{2}
    \lesssim \la^{\frac{1}{2}}\e^{\frac{1}{2}}\lpnorm{e^{-i\la^{-1}tH_V}f_{\la}}{2}{M \times [0,\e^{-1}]} \leq \la^{\frac{1}{2}}\altlpnorm{f_{\la}}{2}
\end{split}
\end{equation}
where in the last inequality we used that $$ \Lpopnorm{e^{-i\la^{-1}tH_V}}{2}{M}{2}{M} = 1.$$
For the second estimate we can use (\ref{Strichartz implies time localized 2}) again, along with (\ref{potential height bound}) to see that
\begin{equation}
\begin{split}
    \mixlpnorm{\beta(-D_t/ \la)w_{\leq \ell}}{2}{q_e}{M \times [0,\e^{-1}]} &\lesssim \la^{-\frac{1}{2}}\e^{-\frac{1}{2}}\lpnorm{F_{\leq \ell}}{2}{M \times [0,\e^{-1}]} = \la^{-\frac{1}{2}}\e^{-\frac{1}{2}}\lpnorm{V_{\leq \ell}w}{2}{M \times [0,\e^{-1}]} \\
    &\leq \ell \la^{-\frac{1}{2}}\e^{-\frac{1}{2}}\lpnorm{\eta(\e t)e^{-i\la^{-1}tH_V}f_{\la}}{2}{M \times [0,\e^{-1}]} \lesssim \ell \la^{-\frac{1}{2}}\e^{-1}\altlpnorm{f_{\la}}{2}
\end{split}
\end{equation}
which is better than what we wanted to show since $\e^{-1} = \log\la$.
Finally, using (\ref{Strichartz implies time localized 1}) we have
\begin{equation} \label{application of inhom Strichartz endpoint to endpoint estimate for large height potential}
    \mixlpnorm{\beta(-D_t/ \la)w_{>\ell}}{2}{q_e}{M \times [0,\e^{-1}]} \lesssim \mixlpnorm{F_{> \ell}}{2}{q_e'}{M \times [0,\e^{-1}]} = \mixlpnorm{V_{> \ell}w}{2}{q_e'}{M \times [0,\e^{-1}]}.
\end{equation}
Note that 
$
    \tfrac{1}{q'_e}= \tfrac{1}{q_e} +\tfrac 2d
$.
 By Holder's inequality,
\begin{equation}
    \lpnorm{V_{> \ell}w(\cdot,t)}{q_e'}{M} \leq \lpnorm{V_{> \ell}}{d/2}{M}\lpnorm{w(\cdot,t)}{q_e}{M} = \delta(\ell)\lpnorm{w(\cdot,t)}{q_e}{M}.
\end{equation}
Plugging this back into (\ref{application of inhom Strichartz endpoint to endpoint estimate for large height potential}) we see that 
\begin{equation} \label{applying Holder after inhom endpoint to endpoint Strichartz for large height potential}
    \mixlpnorm{\beta(-D_t/ \la)w_{>\ell}}{2}{q_e}{M \times [0,\e^{-1}]} \leq  C\mixlpnorm{V_{> \ell}w}{2}{q_e'}{M \times [0,\e^{-1}]} \leq C\delta(\ell)\mixlpnorm{w}{2}{q_e}{M \times [0,\e^{-1}]}.
\end{equation}
Now, choose $\ell$ large enough so that $C \delta(\ell) < \frac{1}{2}$. Then (\ref{time localized endpoint to endpoint estimate for large height potential}) follows immedaitely from (\ref{applying Holder after inhom endpoint to endpoint Strichartz for large height potential}). Thus, the proof of \eqref{rescaled d.1} is complete.

\bibliography{refs}
\bibliographystyle{abbrv}

\end{document}